\documentstyle{amltd}
\begin{document}
\def\joinrel{\mathrel{\mkern-4mu}}
\def\relbar{\mathrel{\smash-}}
\def\lrar{\relbar\joinrel\relbar\joinrel\rightarrow}
\def\xrightarrow#1{\stackrel{#1}{\longrightarrow}}
\def\scrs#1{{\scriptstyle #1}\big\downarrow\phantom{\scriptstyle#1}}
\def\nhs{\hskip-5pt}

\def\dline#1{\stackrel{#1}{=\joinrel=\joinrel=\joinrel=\joinrel=\joinrel=\joinrel=}}
\def\ldline#1{\stackrel{#1}{=\joinrel=\joinrel=\joinrel=\joinrel=\joinrel=\joinrel=\joinrel=\joinrel=\joinrel=\joinrel=\joinrel=\joinrel=\joinrel=\joinrel=\joinrel=\joinrel=}}

 %\input BoxedEps.Tex %for mac
%\SetTexturesEPSFSpecial %for mac
\input boxedeps.tex % unixx
\SetepsfEPSFSpecial % unix
\HideDisplacementBoxes
\def\figin#1#2{
$$
 {\BoxedEPSF{#1.eps scaled
#2}%
}%
$$
\noindent}

\annalsline{158}{2003}
\received{October 23, 2000}
\startingpage{115}
\def\bye{\end{document}}
  \font\tenrm=cmr10
\def\ritem#1{\item[{\rm #1}]}

\def\hensp#1{\quad\hbox{#1}\quad }

\def\Fr{{\rm Fr}}
\def\inv{\inverse}
 %for Bbb in amstex
\catcode`\@=11
\font\twelvemsb=msbm10 scaled 1100
\font\tenmsb=msbm10
%\font\ninemsb=msbm7 scaled 1100%msbm9
\font\ninemsb=msbm10 scaled 800
\newfam\msbfam
\textfont\msbfam=\twelvemsb  \scriptfont\msbfam=\ninemsb
  \scriptscriptfont\msbfam=\ninemsb
\def\msb@{\hexnumber@\msbfam}
\def\Bbb{\relax\ifmmode\let\next\Bbb@\else
 \def\next{\errmessage{Use \string\Bbb\space only in math
mode}}\fi\next}
\def\Bbb@#1{{\Bbb@@{#1}}}
\def\Bbb@@#1{\fam\msbfam#1}
\catcode`\@=12

 \catcode`\@=11
\font\twelveeuf=eufm10 scaled 1100
\font\teneuf=eufm10
\font\nineeuf=eufm7 scaled 1100%eufm9
\newfam\euffam
\textfont\euffam=\twelveeuf  \scriptfont\euffam=\teneuf
  \scriptscriptfont\euffam=\nineeuf
\def\euf@{\hexnumber@\euffam}
\def\frak{\relax\ifmmode\let\next\frak@\else
 \def\next{\errmessage{Use \string\frak\space only in math
mode}}\fi\next}
\def\frak@#1{{\frak@@{#1}}}
\def\frak@@#1{\fam\euffam#1}
\catcode`\@=12

%--------------- Author macros ---------------

 \newcommand{\Ends}{\mathop{\rm Ends}\nolimits}
\def\QI{\mathop{\rm QI}\nolimits}
\def\interior{\mathop{\rm int}\nolimits}
\def\Length{\mathop{\rm Length}\nolimits}
\def\SL{\mathop{\rm SL}\nolimits}
\def\Stab{\mathop{\rm Stab}\nolimits}
\def\closure{\mathop{\rm cl}\nolimits}
\def\St{\mathop{\rm St}\nolimits}
\def\QE{\mathop{\rm QE}\nolimits}
\def\Hull{{\cal H}}
\def\Core{\mathop{\rm Core}\nolimits}
\def\diam{\mathop{\rm Diam}\nolimits}
\def\image{\mathop{\rm image}\nolimits}
\def\midpoint{\mathop{\rm mid}\nolimits}
\def\link{\mathop{\rm Lk}\nolimits}
\def\PD{\mathop{\rm PD}\nolimits}
\def\CR{\rho}
\def\Isom{\mathop{\rm Isom}\nolimits}
\def\Homeo{\mathop{\rm Homeo}\nolimits}
\def\Ker{\mathop{\rm Ker}\nolimits}
\def\kernel{\Ker}
\def\Sim{\mathop{\rm Sim}\nolimits}
\def\height{\mathop{\rm ht}\nolimits}
\def\sep{\mathop{\rm sep}\nolimits}
\def\Vertices{\mathop{\rm Verts}\nolimits}
\def\Edges{\mathop{\rm Edges}\nolimits}
\def\BS{\mathop{\rm BS}\nolimits}
\def\Aut{\mathop{\rm Aut}\nolimits}
\def\cd{\mathop{\rm cd}\nolimits}
\def\supp{\mathop{\rm supp}\nolimits}
\def\gcf{\mathop{\rm gcf}\nolimits}
\def\QC{\mathop{\rm QC}\nolimits}
\def\Conf{\mathop{\rm Conf}\nolimits}
\def\IValence{\mathop{\rm I-Val}\nolimits}
\def\PSL{\mathop{\rm PSL}\nolimits}
\def\Star{\mathop{\rm Star}\nolimits}
\def\GL{\mathop{\rm GL}\nolimits}

\newcommand\R{{\bf R}}
\newcommand\reals{\R}
\newcommand\Q{{\bf Q}}
\newcommand\hyp{{\bf H}}
\newcommand\complex{{\bf C}}
\newcommand\Z{{\bf Z}}
\newcommand\Euc{{\bf E}}
\renewcommand\H{{\cal H}}
\newcommand\F{{\cal F}}
\newcommand\E{{\cal E}}
\newcommand\C{{\cal C}}
\newcommand\K{{\cal K}}
\newcommand\Y{{\cal Y}}

\newcommand\solv{{\scshape solv}}

\newcommand\inject{\hookrightarrow}
\newcommand\Sum{\sum}
\newcommand\infinity{\infty}
\newcommand{\bdy}{\partial}
\newcommand{\from}{\colon}
\newcommand\composed{\circ}
\newcommand\suchthat{\bigm|}
\newcommand\inverse{{-1}}
\newcommand\union{\cup}
\newcommand\Union{\bigcup}
\newcommand\abs[1]{\left| #1 \right|}
\newcommand\card[1]{\abs{#1}}
\newcommand\subgroup{<}
\newcommand\edot{\dot e}
\newcommand\homeo\approx
\newcommand\cbdy\delta

% OLD COARSE MACROS
%\newcommand\cequiv{\approx_c}
%\newcommand\ce\cequiv
%\newcommand\coarsecap{\cap_c}
%\newcommand\cc{\subset_c}
%\newcommand\ccp{\subsetneq_c}
%\newcommand\ceq{=_c}
%\newcommand\coarsely[1]{\,\,[#1]}

% NEW COARSE MACROS
\newcommand\QIClass[1]{\langle\langle#1\rangle\rangle}
\newcommand\csubset[1]{\stackrel{#1}{\subset}}
\newcommand\ceq[1]{\stackrel{#1}{=}}

\newcommand\Id{{{\rm Id}}}
\newcommand\intersect{\cap}

\newcommand\Svarc{\v{S}varc}
\newcommand\restrict{\bigm|}
\newcommand\semidirect{\rtimes}

\newcommand\cross{\times}
\newcommand\into{\hookrightarrow}
\newcommand\Hequiv{\sim}
\newcommand\ext{{\rm int}}
%fix this later
\newcommand\Haus{{\cal H}}
\renewcommand\O{{\rm O}}
\newcommand\M{{\cal M}}
\newcommand\Mobius{Mobi\"us}
\newcommand\Poincare{Poincar\'e}

\newcommand\<\langle
\renewcommand\>\rangle
\newcommand\disjunion{{\textstyle\coprod}}
\newcommand\Disjunion{{\displaystyle\coprod}}
\newcommand\Min{{\min}}
\newcommand\Max{{\max}}

\newcommand\U{{\cal U}}
\newcommand\V{{\cal V}}
\renewcommand\H{{\cal H}}
\newcommand\G{{\cal G}}

\newcommand\hf{\hat{f}}
\newcommand\wt{\widetilde}

%-------------- Author entries --------------------

\title{Quasi-actions on trees I.\\
 Bounded valence}
\shorttitle{Quasi-actions on trees I} 
  \twoauthors{Lee Mosher, Michah Sageev,}{Kevin Whyte}

  \institutions{Rutgers University, Newark, NJ
\\
{\eightpoint {\it E-mail address\/}: mosher@andromeda.rutgers.edu}\\
\vglue6pt
 Technion, Israel Institute  of Technology, Haifa, Israel\\
{\eightpoint {\it E-mail address\/}: sageevm@techunix.technion.ac.il}\\
 \vglue6pt
The University of Chicago, Chicago, IL
\\
{\eightpoint {\it E-mail address\/}: 
kwhyte@math.uchicago.edu}}
%-------------- Article Text--------------------

\hyphenation{quasi-symmetric quasi-conformal quasi-similarity}

 \newcommand{\nb}[1]{#1\nobreakdash-}

\newcommand{\textmatrix}[4]{\bigl( \begin{smallmatrix} #1 & #2 \\ #3 & #4
\end{smallmatrix} \bigr)}

  \centerline{\bf Abstract}
\vglue12pt
Given a bounded valence, bushy tree $T$, we prove that any
cobounded quasi-action of a group $G$ on $T$ is quasiconjugate to an
action of $G$ on another bounded valence, bushy tree $T'$. This theorem
has many applications: quasi-isometric rigidity for fundamental groups of
finite, bushy graphs of coarse $\PD(n)$ groups for each fixed $n$; a
generalization to actions on Cantor sets of Sullivan's theorem about
uniformly quasiconformal actions on the 2-sphere; and a
characterization of locally compact topological groups which contain a
virtually free group as a cocompact lattice. Finally, we give the first
examples of two finitely generated groups which are quasi-isometric and
yet which cannot act on the same proper geodesic metric space, properly
discontinuously and cocompactly by isometries.
 
\vglue-12pt
\section{Introduction}
 
\label{SectionIntro}
 
A {\it quasi\/{\rm -}\/action} of a group $G$ on a metric
space $X$ associates to each $g \in G$ a quasi-isometry $A_g \from x \to
g \cdot x$ of $X$, with uniform quasi-isometry constants, so that
$A_{\Id}=\Id_X$, and so that the distance between $A_g\composed A_h$ and
$A_{gh}$ in the sup norm is uniformly bounded independent of $g,h \in G$.

Quasi-actions arise naturally in geometric group theory: if a metric
space $X$ is quasi-isometric to a finitely generated group $G$ with its
word metric, then the left action of $G$ on itself can be
``quasiconjugated'' to give a quasi-action of $G$ on $X$. Moreover, a
quasi-action which arises in this manner is {\it cobounded} and
{\it proper}; these properties are generalizations of {\it cocompact}
and {\it properly discontinuous} as applied to isometric actions.

Given a metric space $X$, a fundamental problem in geometric group theory
is to characterize groups quasi-isometric to $X$, or equivalently, to
characterize groups which have a proper, cobounded quasi-action on $X$.
A~more general problem is to characterize arbitrary quasi-actions on $X$
up to quasiconjugacy. This problem is completely solved in the
prototypical cases $X=\hyp^2$ or $\hyp^3$: any quasi-action on $\hyp^2$
or $\hyp^3$ is quasiconjugate to an isometric action. When $X$ is an
irreducible symmetric space of nonpositive curvature, or an irreducible
Euclidean building of rank $\ge 2$, then as recounted below similar
results hold, sometimes with restriction to cobounded quasi-actions,
sometimes with stronger conclusions.

The main result of this paper, Theorem~1, gives
a complete solution to the problem for cobounded quasi-actions in the case
when $X$ is a bounded valence tree which is {\it bushy}, meaning coarsely
that the tree is neither a point nor a
line. Theorem~1 says that any cobounded
quasi-action on a bounded valence, bushy tree is quasiconjugate to an
isometric action, on a possibly different tree.

We give various applications of this result.

For instance, while the typical way to prove that two groups are
quasi-isometric is to produce a proper metric space on which they each
have a proper cobounded action, we provide the first examples of two
quasi-isometric groups for which there does not exist {\it any} proper
metric space on which they both act, properly and coboundedly; our
examples are virtually free groups. We do this by determining which
locally compact groups $\G$ can have discrete, cocompact subgroups that
are virtually free of finite rank $\ge 2$: $\G$ is closely related to the
automorphism group of a certain bounded valence, bushy tree $T$. In
\cite{MSW02a} these results are applied to
characterize which trees $T$ are the ``best'' model geometries for
virtually free groups; there is a countable infinity of
``best'' model geometries in an appropriate sense.

Our main application is to quasi-isometric rigidity for homogeneous
graphs of groups; these are finite graphs of finitely generated groups in
which every edge-to-vertex injection has finite index image. For
instance, we prove quasi-isometric rigidity for fundamental groups of
finite graphs of virtual $\Z$'s, and by applying previous results we then
obtain a complete classification of such groups up to quasi-isometry.
More generally, we prove quasi-isometric rigidity for a homogeneous graph
of groups $\Gamma$ whose vertex and edge groups are ``coarse'' $\PD(n)$
groups, as long as the Bass-Serre tree is bushy---any finitely generated
group $H$ quasi-isometric to $\pi_1\Gamma$ is the fundamental group of a
homogeneous graph of groups $\Gamma'$ with bushy Bass-Serre tree whose
vertex and edge groups are quasi-isometric to those of $\Gamma$.

Other applications involve the problem of passing from quasiconformal
boundary actions to conformal actions, where in this case the boundary
is a Cantor set.  
Quasi-actions on $\hyp^3$ are studied via the theorem that any
uniformly quasiconformal action on $S^2=\bdy\hyp^3$ is quasiconformally
conjugate to a conformal action; the countable case of this theorem was
proved by Sullivan \cite{Sul81}, and the general case
by Tukia \cite{Tuk80}.\footnote{The fact that Sullivan's theorem
implies QI-rigidity of $\hyp^3$ was pointed out by Gromov to Sullivan in
the 1980's \cite{Sul}; see also \cite{CC92}.}
Quasi-actions on other hyperbolic symmetric spaces are studied via similar
theorems about uniformly quasiconformal boundary actions, sometimes
requiring that the induced action on the triple space be cocompact, as
recounted below. Using Paulin's formulation of uniform quasiconformality
for the boundary of a Gromov hyperbolic space \cite{Pau96},
we prove that when $B$ is the Cantor set, equipped with a quasiconformal
structure by identifying $B$ with the boundary of a bounded valence bushy
tree, then any uniformly quasiconformal action on $B$ whose induced
action on the triple space is cocompact is quasiconformally conjugate
to a conformal action in the appropriate sense. Unlike the more analytic
proofs for boundaries of rank~1 symmetric spaces, our proofs depend on
the low-dimensional topology methods of
Theorem~1.

Quasi-actions on $\hyp^2$ are studied similarly via the induced actions\break
on $S^1=\bdy\hyp^2$. We are primarily interested in one subcase, a
theorem of Hinkkanen \cite{Hin85} which says that any
uniformly quasi-symmetric group action on $\R = S^1 -
\{\hbox{point}\}$ is quasisymmetrically conjugate to a similarity action
on $\R$; an analogous theorem of Farb and Mosher
\cite{FM99} says that any uniform quasisimilarity group action
on $\R$ is bilipschitz conjugate to a similarity action. We prove a
Cantor set analogue of these results, answering a question posed in
\cite{FM99}: any uniform quasisimilarity action on the
$n$-adic rational numbers $\Q_n$ is bilipschitz conjugate to a similarity
action on some $\Q_m$, with $m$ possibly different from $n$.
 
Theorem~1 has also been applied recently by
A.~Reiter \cite{Rei02} to solve quasi-isometric rigidity
problems for lattices in $p$-adic Lie groups with rank~1 factors, for
instance to show that any finitely generated group quasi-isometric to a
product of bounded valence trees acts on a product of bounded valence
trees.

\demo{Acknowledgements}
 The authors are supported in part by the National Science Foundation: the
first author by NSF grant DMS-9803396; the second author by NSF grant
DMS-989032; and the third author by an NSF Postdoctoral
Research Fellowship.
\enddemo
%\begin{small}
 \centerline{\bf Contents}

\def\sni#1{\vglue1pt\noindent{#1}. }
\def\ssni#1{\vglue1pt\noindent\hskip18pt {#1}.}
\def\pni{\vglue1pt\noindent \phantom{1. } }
\sni{1} Introduction
\vglue1pt\noindent \phantom{1. } Acknowledgements
\sni{2} Statements of results
\ssni{2.1} Theorem 1: Rigidity of quasi-actions on bounded valence, bushy trees
\ssni{2.2} Application: Quasi-isometric rigidity for graphs of coarse\\
\phantom{2.2. \hskip15pt } PD$(n)$ groups
\ssni{2.3}  Application: Actions on Cantor sets
\ssni{2.4} Application: Virtually free, cocompact lattices
\ssni{2.5} Other applications \pagebreak

\noindent {3. } Quasi-edges and the proofs of Theorem 1
\ssni{3.1} Preliminaries 
\ssni{3.2} Setup
\ssni{3.3} Quasi-edges
\ssni{3.4} Construction of the 2-complex $X$
\ssni{3.5} Tracks
\sni{4} Application: Quasi-isometric rigidity for graphs of coarse $\PD(n)$ groups
\ssni{4.1} Bass-Serre theory
\ssni{4.2} Geometrically homogeneous graphs of groups
\ssni{4.3} Weak vertex rigidity
\ssni{4.4} Coarse Poincar\'e duality spaces and groups
\ssni{4.5} Bushy graphs of coarse $\PD(n)$ groups
\sni{5} Application: Actions on Cantor sets
\ssni{5.1} Uniformly quasiconformal actions
\ssni{5.2} Uniform quasisimilarity actions on $n$-adic Cantor sets
\sni{6} Application: Virtually free, cocompact lattices

%\end{small}

\section{Statements of results}

2.1. {\it Theorem {\rm 1:} Rigidity of quasi\/{\rm -}\/actions on bounded valence{\rm ,}
bushy trees}. The simplest nonelementary Gromov hyperbolic metric spaces are
homogeneous simplicial trees $T$ of constant valence $\ge 3$. One novel
feature of such geometries is that there is no best geometric model: all
trees with constant valence $\ge 3$ are quasi-isometric to each other.
Indeed, each such tree is quasi-isometric to any tree $T$ satisfying the
following properties: $T$ has {\it bounded valence}, meaning that
vertices have uniformly finite valence; and $T$ is
{\it bushy}, meaning that each point of $T$ is a uniformly bounded
distance from a vertex having at least 3 unbounded complementary
components. In this paper, each tree $T$ is given a geodesic metric in
which each edge has length~$1$; one effect of this is to identify the
isometry group $\Isom(T)$ with the automorphism group of $T$.

Here is our main theorem:\footnote{Theorem~1 and several of
its applications were first presented in
\cite{MSW00}, which also presents results from
Part 2 of this paper \cite{MSW02b}.}

\nonumproclaim{Theorem 1 {\rm (Rigidity of quasi-actions on bounded valence, bushy trees)}} 
If $G \cross T \to T$ is a cobounded quasi\/{\rm -}\/action of a group $G$ on a
bounded valence{\rm ,} bushy tree $T${\rm ,} then there is a bounded valence{\rm ,} bushy
tree $T'${\rm ,} an isometric action $G \cross T' \to T'${\rm ,} and a
quasiconjugacy $f\from T' \to T$ from the action of $G$ on $T'$ to the
quasi\/{\rm -}\/action of $G$ on $T$.
\endproclaim
 
{\it Remark}. Given quasi-actions of $G$ on metric spaces $X,Y$, a
quasiconjugacy is a quasi-isometry $f \from X \to Y$ which is
{\it coarsely $G$-equivariant} meaning that
$d_{Y}\bigl(f(g\cdot x), g \cdot fx\bigr)$ is uniformly bounded
independent of $g \in G$, $x \in X$. Any coarse inverse for $f$ is also
coarsely $G$-equivariant. We remark that properness and coboundedness are
each invariant under quasiconjugation.
\vglue12pt

Theorem~1 complements similar theorems for
irreducible symmetric spaces and Euclidean buildings. The results for
$\hyp^2$ and $\hyp^3$ were recounted above. When $X = \hyp^n$,
$n \ge 4$ \cite{Tuk86}, and when $X = \complex\hyp^n$, $n
\ge 2$ \cite{Cho96}, every cobounded quasi-action is
quasiconjugate to an action on $X$. Note that if $n \ge 4$ then $\hyp^n$
has a noncobounded quasi-action which is not quasiconjugate to any action on
$\hyp^n$
\cite{Tuk81},
\cite{FS87}. When $X$ is a quaternionic hyperbolic
space or the Cayley hyperbolic plane \cite{Pan89b}, or when $X$ is a
nonpositively curved symmetric space or thick Euclidean building,
irreducible and of rank $\ge 2$ \cite{KL97b}, every
quasi-action is actually a bounded distance from an
action on $X$. Theorem~1 complements the
building result because bounded valence, bushy trees with cocompact
isometry group incorporate thick Euclidean buildings of rank~1. However,
the conclusion of Theorem~1 cannot be as strong
as the results of \cite{Pan89b} and \cite{KL97b}. A
given quasi-action on a bounded valence, bushy tree $T$ may not be
quasiconjugate to an action on the same tree $T$ (see
Corollary~\ref{CorollaryNoGoodModel}); and even if it is, it may not be a
bounded distance from an isometric action on $T$.

The techniques in the proof of
Theorem~1 are quite different from the above
mentioned results. Starting from the induced action of $G$ on $\bdy T$,
first we construct an action on a discrete set, then we attach edges
equivariantly to get an action on a locally finite graph quasiconjugate
to the original quasi-action. This graph need not be a tree, however. We
next attach 2-cells equivariantly to get an action on a locally
finite, simply connected 2-complex quasiconjugate to the original
quasi-action. Finally, using Dunwoody's tracks
\cite{Dun85}, we construct the desired tree action.

Theorem~1 is a very general result, making no
assumptions on properness of the quasi-action, and no
assumptions whatsoever on the group $G$. This freedom facilitates
numerous applications, particularly for improper quasi-actions.

\demo{{\rm 2.2.} Application\/{\rm :}\/ Quasi\/{\rm -}\/isometric rigidity
for graphs of coarse $\PD(n)$ groups}
\vglue6pt
From the proper case of Theorem~1 it
follows that any finitely generated group $G$ quasi-isometric to a free
group is the fundamental group of a finite graph of finite groups, and in
particular $G$ is virtually free; this result is a well-known corollary of
work of Stallings \cite{Sta68} and Dunwoody
\cite{Dun85}. By dropping properness we obtain a much
wider array of quasi-isometric rigidity theorems for certain graphs of
groups.

Let $\Gamma$ be a finite graph of finitely generated groups. There is a
vertex group $\Gamma_v$ for each $v \in \Vertices(\Gamma)$; there is an
edge group $\Gamma_e$ for each $e \in \Edges(e)$; and for each end
$\eta$ of an edge $e$, with $\eta$ incident to the vertex $v(\eta)$,
there is an edge-to-vertex injection $\gamma_\eta\from \Gamma_e \to
\Gamma_{v(\eta)}$. Let $G=\pi_1\Gamma$ be the fundamental group, and let
$G \cross T \to T$ be the action of $G$ on the Bass-Serre tree $T$ of
$\Gamma$. See Section \ref{SectionCoarsePDn} for a brief review of
graphs of groups and Bass-Serre trees.

We say that $\Gamma$ is {\it geometrically homogeneous} if each
edge-to-vertex injection $\gamma_\eta$ has finite index image, or
equivalently $T$ has bounded valence. Other equivalent conditions are
stated in Section \ref{SectionCoarsePDn}.

Consider for example the class of \Poincare\ duality $n$ groups or
$\PD(n)$ groups. If $n$ is fixed then any finite graph of virtual $\PD(n)$
groups is geometrically homogeneous, because a subgroup of a $\PD(n)$
group $K$ is itself $\PD(n)$ if and only it has finite index in
$K$ \cite{Bro82}. In particular, if each vertex and edge group
of
$\Gamma$ is the fundamental group of a closed, aspherical manifold of
constant dimension~$n$ then
$\Gamma$ is geometrically homogeneous.

Our main result, Theorem~\ref{TheoremCoarsePDnGraphs}, is stated in terms
of the (presumably) more general class of ``coarse $\PD(n)$ groups''
defined in Section \ref{SectionCoarsePDn}---such groups respond well to
analysis using methods of coarse algebraic topology introduced in
\cite{FS96} and further developed in
\cite{KK99}.  Coarse PD$(n)$ groups include fundamental groups of compact, aspherical manifolds, groups which are
virtually ${\rm PD}(n)$ of  finite type, and all of Davis' examples in [Dav98].  The definition of coarse ${\rm PD}(n)$
being somewhat technical, we defer the definition to Section~4.4.
\enddemo

\proclaimtitle{QI-rigidity for graphs of coarse PD({\it n})
groups}
\specialnumber{2}
\proclaim{Theorem} 
\label{TheoremCoarsePDnGraphs} Given $n\ge 0${\rm ,} if $\Gamma$ is a finite graph of groups
with bushy Bass\/{\rm -}\/Serre tree{\rm ,} such that each vertex and edge group is
  a coarse $\PD(n)$ group{\rm ,} and if $G$ is a finitely
generated group quasi\/{\rm -}\/isometric to $\pi_1\Gamma${\rm ,} then $G$ is the
fundamental group of a graph of groups with bushy Bass\/{\rm -}\/Serre tree{\rm ,} and
with vertex and edge groups quasi\/{\rm -}\/isometric to those of $\Gamma$.
\endproclaim

Another proof of this result was found, later and independently, by
P.\ Papasoglu \cite{Pap02}.

Given a homogeneous graph of groups $\Gamma$, the Bass-Serre tree $T$
satisfies a trichotomy: it is either finite, quasi-isometric to a line,
or bushy \cite{BK90}. Once $\Gamma$ has been reduced so as to have
no valence~1 vertex with an index~1 edge-to-vertex injection, then: $T$
is finite if and only if it is a point, which happens if and only if
$\Gamma$ is a point; and $T$ is quasi-isometric to a line if and only if
it is a line, which happens if and only if $\Gamma$ is a circle with
isomorphic edge-to-vertex injections all around or an arc with isomorphic
edge-to-vertex injections at any vertex in the interior of the arc and
index~2 injections at the endpoints of the arc. Thus, in some sense
bushiness of the Bass-Serre tree is generic.

Theorem \ref{TheoremCoarsePDnGraphs} suggests the following problem. Given
$\Gamma$ as in Theorem~\ref{TheoremCoarsePDnGraphs}, all edge-to-vertex
injections are quasi-isometries. Given $\C$, a quasi-isometry class of
coarse $\PD(n)$ groups, let $\Gamma\C$ be the class of fundamental groups
of finite graphs of groups with vertex and edge groups in $\C$ and with
bushy Bass-Serre tree. Theorem \ref{TheoremCoarsePDnGraphs} says that
$\Gamma\C$ is closed up to quasi-isometry.
 
\vglue6pt {\elevensc Problem 3.} {\it Given $\C${\rm ,} describe the quasi\/{\rm -}\/isometry classes within
$\Gamma\C$.}

\vglue6pt

Here is  a rundown of the cases for which the solution to this
problem is known to us. Given a metric space $X$, such as a finitely
generated group with the word metric, let $\QIClass{X}$ denote the
class of finitely generated groups quasi-isometric to $X$.

Coarse $\PD(0)$ groups are finite groups, and in this case
Theorem~\ref{TheoremCoarsePDnGraphs} reduces to the fact that
$$\QIClass{F_n} = \Gamma\{\hbox{finite groups}\} = \{\hbox{virtual $F_n$
groups, $n \ge 2$}\}
$$
where the notation $F_n$ will always mean the free group of rank $n \ge
2$.

Coarse $\PD(1)$ groups form a single quasi-isometry class $\C =
\QIClass{\Z} = \{\hbox{virtual $\Z$ groups}\}$. By combining work of
Farb and Mosher
\cite{FM98}, \cite{FM99} with work of Whyte
\cite{Why02}, the groups in $\Gamma\QIClass{\Z}$ are classified as
follows:    

\specialnumber{4}\nonumproclaim{Theorem 4 {\rm (Graphs of coarse  {\rm PD}(1)  groups)}}   
If the finitely generated group $G$ is quasi\/{\rm -}\/isometric to a finite graph
of virtual $\Z$\/{\rm '}\/s with bushy Bass\/{\rm -}\/Serre tree{\rm ,} then exactly one of the
following happens\/{\rm :}\/
\begin{itemize}
\item There exists a unique power free integer $n \ge 2$ such that $G$
modulo some finite normal subgroup is abstractly commensurable to the
solvable\break Baumslag\/{\rm -}\/Solitar group $\BS(1,n) = \< a,t \suchthat tat^\inv =
a^n\>$.
\item $G$ is quasi\/{\rm -}\/isometric to any of the nonsolvable
Baumslag\/{\rm -}\/Solitar groups $\BS(m,n) = \<a,t\suchthat ta^mt^\inv=a^n\>$ with
$2\le m < n$.
\item $G$ is quasi\/{\rm -}\/isometric to any group $F \cross \Z$ where $F$ is free
of finite rank $\ge 2$.
\end{itemize}

\endproclaim

\demo{Proof}
By Theorem~\ref{TheoremCoarsePDnGraphs} we have $G=\pi_1\Gamma$ where
$\Gamma$ is a finite graph of virtual $\Z$'s with bushy Bass-Serre tree.
If $G$ is amenable then the first alternative holds, by
\cite{FM99}. If
$G$ is nonamenable then either the second or the third alternative holds,
by \cite{Why02}.
\enddemo
 
For $\C = \QIClass{\Z^n}$, the amenable groups in $\Gamma\QIClass{\Z^n}$
form a quasi-isometrically closed subclass which is classified
up to quasi-isometry in \cite{FM00}, as follows.
By applying Theorem~1 it is shown that
each such group is virtually an ascending HNN group of the form $\Z^n *_M$
where $M
\in \GL(n,\R)$ has integer entries and $\abs{\det(M)} \ge 2$; the
classification theorem of \cite{FM00} says that the absolute
Jordan form of $M$, up to an integer power, is a complete quasi-isometry
invariant. For general groups in $\Gamma\QIClass{\Z^n}$, Whyte reduces the
problem to understanding when two subgroups of
$\GL_n(\reals)$ are at finite Hausdorff distance
\cite{Why}.

For $\C = \QIClass{\hyp^2}$, the subclass of $\Gamma\QIClass{\hyp^2}$
consisting of word hyperbolic surface-by-free groups is
quasi-isometrically rigid and is classified by Farb and Mosher in
\cite{FM02}. The broader classification in
$\Gamma\QIClass{\hyp^2}$ is open.

If $\C$ is the quasi-isometry class of cocompact lattices in an
irreducible, semisimple Lie group $L$ with finite center,
$L \ne\PSL(2,\reals)$, then combining Mostow Rigidity for $L$ with
quasi-isometric rigidity (see
\cite{Far97} for a survey) it follows that for each $G \in \C$
there exists a homomorphism $G \to L$ with finite kernel and discrete,
cocompact image, and this homomorphism is unique up to post-composition
with an inner automorphism of $L$. Combining this with Theorem
\ref{TheoremCoarsePDnGraphs} it follows that $\Gamma\C$ is a single
quasi-isometry class, represented by the cartesian product of any group
in $\C$ with any free group of rank $\ge 2$.

\demo{{R}emark} In \cite{FM99} it is proved that any
finitely generated group $G$ quasi-isometric to $\BS(1,n)$, where $n \ge
2$ is a power free integer, has a finite subgroup $F$ so that $G/F$ is
abstractly commensurable to $\BS(1,n)$.
Theorem~4 can be applied to give a (mostly)
new proof, whose details are found in \cite{FM00}.
\enddemo

2.3. {\it Application\/{\rm :} Actions on Cantor sets}.

\vglue6pt {\it Quasiconformal actions}. The boundary
of a $\delta$-hyperbolic metric space $X$ carries a quasiconformal
structure and a well-behaved notion of uniformly quasiconformal
homeomorphisms, which as Paulin showed can be characterized in terms of
cross ratios \cite{Pau96}; we review this in Section
\ref{SectionCantorActions}. As such, one ask can for a generalization of
the Sullivan-Tukia theorem for $\hyp^3$: is every uniformly
quasiconformal group action on $\bdy X$ quasiconformally conjugate to a
conformal action?

A bounded valence, bushy tree $T$ has Gromov boundary $B=\bdy T$
homeomorphic to a Cantor set, and for actions with an appropriate
cocompactness property we answer the above question in the affirmative for
$B$, where ``conformal action'' is interpreted as the induced action at
infinity of an isometric action on some other bounded valence, bushy
tree. Recall that an isometric group action on a $\delta$-hyperbolic
metric space $X$ is cocompact if and only if the induced action on the
space of distinct triples in $\bdy X$ is cocompact, and the action on $X$
has bounded orbits if and only if the induced action on the space of
distinct pairs in $\bdy X$ has precompact orbits.

\proclaimtitle{Quasiconformal actions on Cantor sets}
\specialnumber{5}
\proclaim{Theorem} 
\label{TheoremQC} If the Cantor set $B$ is equipped with a
quasiconformal structure by identifying $B = \bdy T$ for some bounded
valence{\rm ,} bushy tree $T${\rm ,} if $G \cross B \to B$ is a uniformly
quasiconformal action of a group $G$ on $B${\rm ,} and if the action of $G$ on
the triple space of $B$ is cocompact{\rm ,} then there exists a tree
$T'$ and a quasiconformal homeomorphism $\phi \from B \to \bdy T'$ which
conjugates the $G$\/{\rm -}\/action on $B$ to an action on $\bdy T'$ which is
induced by some cocompact{\rm ,} isometric action of $G$ on $T'$.
\endproclaim
 
\specialnumber{6}\proclaim{{C}orollary}
\label{CorollaryQC}
Under the same hypotheses as Theorem~{\rm \ref{TheoremQC},} $G$ is the
fundamental group of a finite graph of groups $\Gamma$ with finite index
edge\/{\rm -}\/to\/{\rm -}\/vertex injections\/{\rm ;} moreover a subgroup $H \subgroup G$ stabilizes
some vertex of the Bass\/{\rm -}\/Serre tree of $\Gamma$ if and only if the action
of $H$ on the space of distinct pairs in $B$ has precompact orbits.
\endproclaim

Once the definitions are reviewed, the proofs of Theorem
\ref{TheoremQC} and Corollary~\ref{CorollaryQC} are very quick
applications of Theorem~1.

Theorem~\ref{TheoremQC} complements similar theorems for the boundaries of
all rank~1 symmetric spaces. Any uniformly quasiconformal action on the
boundary of $\hyp^2$ or $\hyp^3$ is quasiconformally conjugate to a
conformal action. Any uniformly quasiconformal action on the boundary of
$\hyp^n$, $n \ge 4$ \cite{Tuk86} or of $\complex\hyp^n$
\cite{Cho96}, such that the induced action on the triple space
of the boundary is cobounded, is quasiconformally conjugate to a conformal
action. Any quasiconformal map on the boundary of a quaternionic
hyperbolic space or the Cayley hyperbolic plane is conformal
\cite{Pan89b}.

Also, convergence actions of groups on Cantor sets have been studied in
unpublished work of Gerasimov and in work of Bowditch
\cite{Bow02}. These works show that if the group $G$
has a minimal convergence action on a Cantor set $C$, and if $G$ satisfies
some mild finiteness hypotheses, then there is a $G$-equivariant
homeomorphism between $C$ and the space of ends of $G$.
Theorem~\ref{TheoremQC} and the corollary are in the same vein, though
for a different class of actions on Cantor sets.

\vglue6pt {\it Uniform quasisimilarity actions on the $n$\/{\rm -}\/adics}.  Given $n \ge
2$, let $\Q_n$ be the $n$-adics, a complete metric space whose points are
formal series 
$$\xi=\sum_{i=k}^{+\infinity} \xi_i \, n^i, \quad\hbox{where}\quad
\xi_i\in\Z/n\Z \quad\hbox{and}\quad k \in \Z.
$$
The distance between $\xi,\eta \in \Q_n$ equals $n^{-I}$ where $I$ is the
greatest element of $\Z \union \{+\infinity\}$ such that $\xi_i=\eta_i$
for all $i \le I$. The metric space $\Q_n$ has Hausdorff dimension $1$,
and it is homeomorphic to a Cantor set minus a point.

Given integers $m,n \ge 2$, Cooper proved that the metric spaces $\Q_m,
\Q_n$ are bilipschitz equivalent if and only if there exists integers $k
\ge 2$, $i,j \ge 1$ such that $m=k^i$, $n=k^j$ (see Cooper's appendix
to \cite{FM98}). Thus, each bilipschitz class of
$n$-adic metric spaces is represented uniquely by some $\Q_m$ where $m$ is
not a proper power.

A {\it similarity} of a metric space $X$ is a bijection $f \from X
\to X$ such that the ratio $d(f\xi,f\eta)/d(\xi,\eta)$ is
constant, over all $\xi \ne \eta \in X$. A
$K$-quasisimilarity, $K \ge 1$, is a bijection
$f\from X \to X$ with the property that
$$
\frac{d(f\zeta,f\omega)}{d(\zeta,\omega)} \Biggm/
\frac{d(f\xi,f\eta)}{d(\xi,\eta)}
\le K, \quad\hbox{for all } \zeta \ne \omega, \xi \ne \eta \in X.
$$
A $1$-quasisimilarity is the same thing as a similarity.

In \cite{FM99} it was asked whether any uniform
quasisimilarity action on $\Q_n$ is bilipschitz conjugate to a similarity
action, as long as $n$ is not a proper power. In retrospect this is not
quite the correct question, and in fact there is an easy counterexample:
the full similarity group of $\Q_4$ acts as a uniform quasisimilarity
group on $\Q_2$, but there is no bilipschitz conjugacy to a similarity
action on $\Q_2$ (see the end of \S5.2 for details). This can
be extended to show that for any $m$ and any $i \ne j$, there is a uniform
quasisimilarity action on $\Q_{m^i}$ which is not bilipschitz conjugate
to a similarity action on $\Q_{m^j}$. The following theorem resolves the
issue in the best possible way, at least for actions satisfying the
appropriate cocompactness property, which for a ``punctured'' Cantor set
means cocompactness on the set of distinct pairs:
 
\vglue4pt
{\elevensc Theorem 7} (Uniform quasisimilarity actions on $n$-adic Cantor sets).
{\it Given $n \ge 2${\rm ,} suppose that $G \cross\Q_n \mapsto \Q_n$ is a uniform
quasisimilarity action\/{\rm :} there exists $K \ge 1$ such that each element of
$G$ acts by a $K$\/{\rm -}\/quasisimilarity. Suppose in addition that the induced
action of $G$ on the space of distinct pairs in $\Q_n$ is
cocompact. Then there exists $m \ge 2$ and a bilipschitz homeomorphism
$\Q_n \mapsto \Q_m$ which conjugates the $G$ action on $\Q_n$ to a
similarity action on $\Q_m$.}
\vglue4pt 

This theorem generalizes the similar result of \cite{FM99} for
uniform quasisimilarity actions on $\R$, which was in turn an analogue of
Hinkkanen's theorem \cite{Hin85} for uniformly
quasisymmetric actions on $\R$. We do not know whether the cocompactness
hypothesis is necessary, but it is a useful and commonly occurring
boundedness property.

As with Theorem \ref{TheoremQC}, the result of Theorem 7
allows us to make some extra conclusions about the algebraic structure of
the given group $G$, namely an ascending HNN structure whose base group
is geometrically constrained:
 \vglue3pt
{\elevensc Corollary 8.}
 {\it Under the same hypotheses as Theorem~{\rm 7,}
there is an ascending {\rm HNN} decomposition
$$G = \<H,t \suchthat tht^\inv = \phi(h), \forall h \in H\>,
$$ 
where $H$ is a subgroup of $G$ and $t \in G${\rm ,} such that $\phi \from H \to
H$ is a self\/{\rm -}\/monomorphism with finite index image{\rm ,} and the action of $H$
on $\Q_n$ is uniformly bilipschitz.}
 
\vglue4pt
2.4. {\it Application\/{\rm :} Virtually free{\rm ,} cocompact lattices}.
\vglue2pt
Given a finitely generated group $G$, one can ask to describe the
{\it model geometries} for $G$, the proper metric spaces $X$ on which
$G$ acts, properly and coboundedly by isometries. More generally,
motivated by $\Isom(X)$, one can ask to describe the locally compact
topological groups $\Gamma$ for which there is a discrete, cocompact,
virtually faithful representation $G \to\Gamma$. Then, given a
quasi-isometry class $\C$ of finitely generated groups, one can ask:
\begin{itemize}
\item Is there a common model geometry $X$ for every group in $\C$?
\item Is there a common locally compact group $\Gamma$, in which every
group of $\C$ has a discrete, cocompact, virtually faithful
representation?
\end{itemize}

For example, the Sullivan-Tukia theorem answers these two questions
affirmatively for the quasi-isometry class $\QIClass{\hyp^n}$, using the
space $X=\hyp^n$ and the group $\Gamma=\Isom(\hyp^n)$. Most
quasi-isometric rigidity theorems in the literature provide similarly
affirmative answers for the quasi-isometry class under consideration,
e.g.\ \cite{BP00}, \cite{FM99},
\cite{FM00}, \cite{FM02}, \cite{KL97a},
\cite{KL97b}, \cite{Tab00}; see \cite{Far97}
for a survey. On the other hand, it may be true that there are no
less than two model geometries for the quasi-isometry class
$\QIClass{\hyp^2\cross\reals} = \QIClass{\widetilde{\PSL_2\reals}}$
\cite{Rie01}.

We show that the above questions have a negative answer for the
quasi-isometry class $\QIClass{F_n}$. Our main tool is the following:

\proclaimtitle{Virtually free, cocompact lattices}
\specialnumber{9}\proclaim{Theorem} 
\label{TheoremLattices}
Let $\G$ be a locally compact topological group which contains a
cocompact lattice in the class $\QIClass{F_n}$. Then there exists a
cocompact action of $\G$ on a bushy tree $T$ of bounded valence{\rm ,} inducing
a continuous{\rm ,} proper homomorphism $\G \to \Isom(T)$ with compact kernel
and cocompact image.
\endproclaim

\specialnumber{10}\proclaim{{C}orollary}
\label{CorollaryNoGoodModel}
There exist groups $G,G' \in \QIClass{F_n}$ such that\/{\rm :}\/
\begin{itemize}
\item $G,G'$ do not act properly discontinuously and cocompactly by
isometries on the same proper geodesic metric space.
\item $G,G'$ do not have discrete{\rm ,} cocompact{\rm ,} virtually faithful
representations into the same locally compact group.
\end{itemize}

\endproclaim

Theorem \ref{TheoremLattices} reduces the corollary to the statement
that there are virtually free groups which cannot act properly and
cocompactly on the same tree.  It is easy to produce examples of this
phenomenon, for example $\Z / p * \Z / p$ and $\Z / q * \Z  / q$, for
distinct primes $p,q \ge 3$. These are the first examples of
quasi-isometric groups which are
known not to have a common geometric model.

Theorem~\ref{TheoremLattices} complements a recent result of Alex Furman
\cite{Fur01} concerning an irreducible lattice $G$ in a
semisimple Lie group~$\Gamma$. Furman's result shows, {\it except} when
$G \in \QIClass{F_n}$ is noncocompact in $\SL(2,\R)$, that any locally
compact group $\G$ in which $G$ has a discrete, cocompact, virtually
faithful representation is very closely related to the given Lie group
$\Gamma$.

{\it Remark}. The techniques of the above results should
apply to more general homogeneous graphs of groups. In particular, one
ought to be able to  determine, using these ideas, which of the
Baumslag-Solitar groups $\BS(m,n)$ (\cite{Why02}) are cocompact
lattices in the same locally compact group. Also, using the computation
of $\QI(\BS(1,n))$ in \cite{FM98}, it should be possible to
give a conjugacy classification of the maximal uniform cobounded
subgroups of $\QI(\BS(1,n))$, analogous to
Theorem~\ref{TheoremModelGeometries} below.
 \vglue12pt

2.5. {\it Other applications}.

\vglue8pt {\it Quasi\/{\rm -}\/actions on products of trees}. Recently A.~Reiter
\cite{Rei02} has combined
Theorem~1 with results of Kleiner and Leeb on
quasi-isometric rigidity for Euclidean buildings
\cite{KL97b} to prove:

\specialnumber{11}\proclaim{Theorem}
Suppose that $G$ is a finitely generated group quasi\/{\rm -}\/isometric to a
product of trees $\Pi_{i=1}^k T_i${\rm ,} each tree of bounded valence. Then
$G$ has a finite index subgroup of index at most $k!$ which is isomorphic
to a discrete cocompact subgroup of $\Isom\left(\Pi_{i=1}^k T'_i\right)$
where each $T'_i$ is a tree quasi\/{\rm -}\/isometric to the
corresponding~$T_i$. \endproclaim

The finite index arising in this theorem comes from the fact that $G$ is
allowed to permute the $k$ factors among themselves. For example, every
group quasi-isometric to a product of two bounded valence bushy trees has
a subgroup of index $\le 2$ which acts, properly and coboundedly, on a
product of two bounded valence, bushy trees. This quasi-isometry class
contains all products of two free groups of rank $\ge 2$, but it also
contains torsion free simple groups \cite{BM97}.

\vglue12pt {\it Maximally symmetric trees.}
In light of Theorem~\ref{CorollaryNoGoodModel} showing that there is no
single model geometry for an entire quasi-isometry class $\QIClass{F_n}$,
one might still ask for a list of the ``best'' model geometries for the
class. In \cite{MSW02a} we apply
Theorem~1 to show that these consist of certain
trees which are ``maximally symmetric''.

Recall that for any metric space $X$ the {\it quasi-isometry group}
$\QI(X)$ is the group of self quasi-isometries of $X$ modulo
identification of quasi-isometries which have bounded distance in the sup
norm. A subgroup $H \subgroup \QI(X)$ is {\it uniform} if it can be
represented by a quasi-action on $X$. A uniform subgroup $H \subgroup
\QI(X)$ is {\it cobounded} if the induced quasi-action of $H$ on $X$ is
cobounded. A bounded valence, bushy tree $T$ is
{\it cocompact} if $\Isom(T)$ acts cocompactly on $T$; equivalently, the
image of the natural homomorphism $\Isom(T) \to \QI(T)$ is cobounded. We
say that $T$ is {\it minimal} if it has no valence~1 vertices; minimality
implies that $\Isom(T)$ has no compact normal subgroups, and that
$\Isom(T)\to \QI(T)$ is injective, among other nice properties.

\proclaimtitle{Characterizing maximally symmetric trees
\cite{MSW02a}}
\specialnumber{12}\proclaim{Theorem} 
\label{TheoremMaximal}
For any bounded valence{\rm ,} bushy{\rm ,} cocompact{\rm ,} minimal tree $T${\rm ,} the
following are equivalent\/{\rm :}
\vglue2pt
\noindent\hskip16pt\hangindent=31pt\hangafter=1 $\bullet$ \enspace
 $\Isom(T)$ is a maximal uniform cobounded subgroup of $\QI(T)$.
\vglue2pt
\noindent\hskip16pt\hangindent=31pt\hangafter=1 $\bullet$ \enspace For any bounded valence{\rm ,} bushy{\rm ,}
minimal tree $T'${\rm ,} any continuous{\rm ,} proper{\rm ,} cocompact embedding $\Isom(T) \to \Isom(T')$ is
an isomorphism.
\vglue2pt

\noindent\hskip16pt\hangindent=31pt\hangafter=1 $\bullet$ \enspace
For any locally compact group $\G$ without compact normal subgroups{\rm ,} any
continuous{\rm ,} proper{\rm ,} cocompact embedding $\Isom(T) \to \G$ is an
isomorphism. 

\endproclaim

Such trees $T$ are called {\it maximally symmetric}.
Theorem~\ref{TheoremMaximal} says nothing about existence of maximally
symmetric trees. In \cite{MSW02a} we also
prove:

 \proclaimtitle{Enumerating maximally symmetric trees}
\specialnumber{13}\proclaim{Theorem}
\label{TheoremModelGeometries}
Fix a bounded valence{\rm ,} bushy tree $\tau$. Every uniform
cobounded subgroup of $\QI(\tau)$ is contained in a maximal uniform
cobounded subgroup. Every maximal uniform cobounded subgroup of
$\QI(\tau)$ is identified with the isometry group of some maximally
symmetric tree $T$ via a quasi\/{\rm -}\/isometry $T\leftrightarrow\tau${\rm ,}
inducing a natural one\/{\rm -}\/to\/{\rm -}\/one correspondence between conjugacy classes of
maximal uniform cobounded subgroups of $\QI(\tau)$ and isometry classes of
reduced maximally symmetric trees $T$. There is a countable infinity of
such isometry classes\/{\rm ;} and there is a countable infinity of these isometry
classes represented by trees $T$ which support a proper{\rm ,} cobounded
group action.
\endproclaim

The term ``reduced'' refers to a simple combinatorial operation that
simplifies maximally symmetric trees, as explained in
\cite{MSW02a}.

To summarize, there is a countable infinity of ``best'' geometries for
the quasi-isometry class $\QIClass{F_n}$, distinct up to isometry.
Examples include any homogeneous tree of constant valence $\ge 3$, and
any bipartite, bihomogeneous tree of valences $p\ne q \ge 3$; these each
have proper, cobounded actions, and there is still a countable infinity
of other examples. 

Theorem~\ref{TheoremModelGeometries} should be contrasted with the fact
that if $X$ is a rank~1 symmetric space then $\QI(X)$ has
a unique maximal uniform cobounded subgroup up to conjugacy, namely
$\Isom(X)$; this follows from the fact that every cobounded, uniformly
quasiconformal subgroup acting on $\bdy X$ is quasiconformally conjugate
to a conformal group.

\vglue-12pt

\section{Quasi-edges and the proof of Theorem~1}
\label{SectionQuasiedges}
%\input{QTPartI_Quasiedges.tex}
% QTPartI_Quasiedges.tex
\vglue-6pt

3.1. {\it Preliminaries}.

\vglue4pt {\it Coarse language}. Let $X$ be a metric space. Given $A \subset X$
and $R \ge 0$, denote $N_R(A) = \{x \in X \suchthat \exists a \in A
\quad\hbox{such that}\quad d(a,x) \le R\}$. Given subsets $A, B \subset
X$, let $A \csubset{[R]} B$ denote $A \subset N_R(B)$. Let
$A\csubset{c} B$ denote the existence of $R \ge 0$ such that $A
\csubset{[R]} B$; this is called {\it coarse containment} of $A$ in $B$.
Let $A \ceq{[R]} B$ denote the conjunction of $A \csubset{[R]} B$ and
$B\csubset{[R]} A$; this is equivalent to the statement $d_\Haus(A,B) \le
R$ where $d_\Haus(\cdot,\cdot)$ denotes Hausdorff distance. Let $A \ceq{c}
B$ denote the existence of $R$ such that
$A\ceq{[R]} B$; this is called {\it coarse equivalence} of
$A$ and $B$. 

Given a metric space $X$ and subsets $A,B$, we say that a subset $C$ is a
{\it coarse intersection} of $A$ and $B$ if we have $N_R(A) \intersect
N_R(B) \ceq{c} C$ for all sufficiently large
$R$. A coarse intersection of $A$ and $B$ may not exist, but
if one does exist then it is well-defined up to coarse equivalence.

Given metric spaces $X,Y$, a map $f \from X \to Y$ is {\it coarse
Lipschitz} if $f$ stretches distances by at most an affine function:
there exist $K \ge 1$, $C \ge 0$ such that
$$d_Y(fx,fy) \le K d_X(x,y) + C
$$
We say that $f$ is a {\it uniformly proper embedding} if, in addition,
$f$ compresses distances by a uniform amount: there exists a proper,
increasing function $\rho\from [0,\infinity) \to [0,\infinity)$ such that
$$\rho(d_X(x,y)) \le d_Y(fx,fy)
$$
More precisely we say that $f$ is a $(K,C,\rho)$-uniformly proper
embedding. If we can take $\rho(d)=\frac{1}{K} d - C$ then we say that
$f$ is a {\it $K,C$ quasi-isometric embedding}. If furthermore
$f(X)\ceq{[C]} Y$ then we say that $f$ is a {\it $K,C$ quasi-isometry
between $X$ and $Y$}. A {\it $C'$-coarse inverse} of $f$ is a $K,C'$
quasi-isometry $g \from Y\to X$ such that $x
\ceq{[C']} g(f(x))$ and $y \ceq{[C']} f(g(y))$, for all $x \in X$, $y\in
Y$. A simple fact says that for all $K,C$ there exists $C'$ such that
each $K,C$ quasi-isometry has a $C'$-coarse inverse.

Let $G$ be a group and $X$ a metric space. A {\it $K,C$ quasi-action} of
$G$ on $X$ is a map $G \cross X \to X$, denoted $(g,x)
\mapsto A_g(x)=g\cdot x$, so that for each $g\in G$ the map
$A_g \from X \to X$ is a $K,C$ quasi-isometry of $X$, and for each $x \in
X, g,h \in G$ we have
$$g \cdot (h \cdot x) \ceq{[C]} (gh) \cdot x 
$$
In other words, the sup norm distance between $A_g \composed A_h$ and
$A_{gh}$ is at most $C$. A quasi-action is {\it cobounded} if there
exists a constant $R$ such that for each $x \in X$ we have $G \cdot x
\ceq{[R]} X$. A quasi-action is {\it proper} if for each $R$
there exists $M$ such that for all $x,y\in X$, the cardinality of the set
$$\{g
\in G \suchthat \bigl(g \cdot N(x,R)\bigr) \intersect N(y,R) \ne
\emptyset\}$$ is at most $M$. Note that if $G \cross X \to X$ is an
isometric action on a proper metric space, then ``cobounded'' is
equivalent to ``cocompact'' and ``proper'' is equivalent to ``properly
discontinuous''.

Given a group $G$ and quasi-actions of $G$ on metric spaces $X,Y$, a
{\it quasiconjugacy} is a quasi-isometry $f \from X \to Y$ such that for
some $C \ge 0$ we have $f(g \cdot x) \ceq{[C]} g \cdot fx$ for all $g \in
G$, $x  \in X$. Properness and coboundedness are invariants of
quasiconjugacy.

A fundamental principle of geometric group theory says that if $G$ is a
finitely generated group equipped with the word metric, and if $X$ is a
proper geodesic metric space on which $G$ acts properly discontinuously
and cocompactly by isometries, then $G$ is quasi-isometric to $X$.

A partial converse to this result is the {\it quasi-action principle}
which says that if $G$ is a finitely generated group with the word metric
and $X$ is a metric space quasi-isometric to $G$ then there is a
cobounded, proper quasi-action of $G$ on $X$; the constants for this
quasi-action depend only on the quasi-isometry constants between $G$ and
$X$.

\vfil {\it Ends.} Recall the {\it end compactification} of a locally
compact space\break Hausdorff~$X$. The direct system of compact subsets of $X$
under inclusion has a corresponding inverse system of unbounded
complementary components of compact sets, and an {\it end} is an element
of the inverse limit. Letting $\Ends(X)$ be the set of ends, there is a
compact Hausdorff topology on $\overline X = X \union \Ends(X)$ in which
$X$ forms a dense open set, where for each $e \in \Ends(X)$ there is one
basic open neighborhood of $e$ for each compact subset $K \subset X$,
consisting of the unbounded component $U$ of $X-K$ corresponding to $e$
together with all ends $e'$ for which $U$ is the corresponding
unbounded component of $X-K$.

If $f \from X \to Y$ is a quasi-isometry between proper geodesic metric
spaces then there is a natural induced homeomorphism $\Ends(X) \to
\Ends(Y)$.

If $T$ is a bounded valence, bushy tree then $\Ends(T)$ is a Cantor set.
Moreover, there is a natural homeomorphism between $\Ends(T)$ and the
Gromov boundary of $T$.

\vfil 3.2. {\it Setup}.
\vfil
Let $G$ be a finitely generated group quasi-acting on a bounded valence,
bushy tree $T$. Assume that the quasi-action is cobounded. To prove
Theorem~1 we must construct a quasiconjugacy to
an isometric action of $G$ on another tree.

Before continuing, we immediately reduce to the case where every vertex of
$T$ has valence $\ge 3$. To do this we need only construct a
quasi-isometry $\phi$ from $T$ to a bounded valence tree in which every
vertex  has valence $\ge 3$, for we can then use $\phi$ to quasiconjugate
the given $G$-quasi-action on $T$.

Let $\beta$ be a bushiness constant for $T$: every vertex of $T$ is
within distance $\beta$ of a vertex with $\ge 3$ complementary components.
There is a $\beta$-bushy subtree $T' \subset T$ containing \pagebreak  no valence
$1$ vertices, such that every vertex of $T$ is within distance $\beta$ of
a vertex of $T'$. The nearest point projection map $T \to T'$ is a
$(1,\beta)$ quasi-isometry. Since each valence $2$ vertex of $T'$ is
within distance $\beta$ of a vertex of valence $\ge 3$, we may next
produce a tree $T''$ by changing the tree structure on $T'$, removing
vertices of valence $2$ and conglomerating any path of edges through
valence $2$ vertices into a single edge. The ``identity'' map $T' \to
T''$ is a $\beta,\beta$ quasi-isometry, and $T''$ is the desired tree
in which each vertex has valence $\ge 3$. Replacing $T$ by $T''$, we may
henceforth assume every vertex of $T$ has valence $\ge 3$.

While our ultimate goal is a quasiconjugacy to an action on a tree, our
intermediate goal will be a quasiconjugacy to an action on a certain
2-complex: we construct an isometric action of $G$ on a
2-complex $X$, and a quasiconjugacy $f \from X \to T$, so that $X$ is
simply connected and uniformly locally finite. Once this is accomplished
we use Dunwoody tracks to construct a quasiconjugacy from the $G$ action
on $X$ to a $G$ action on a tree.

%with the
%following properties:
%\begin{itemize}
%\item $X$ is simply-connected.
%\item $X$ is uniformly locally finite.
%\item $f \from X \to T$ is {\it tight}, meaning that $f(X^0) \subset
%T^0$, $f$ maps each edge of $X$ to either a vertex of $T$ or a geodesic
%in $T$, and $f$ maps each 2-simplex of $X$ in the manner described
%in Figure~\ref{FigureSimplexMap}.
%\end{itemize}

%\begin{figure}
%\centeredepsfbox{SimplexMap.eps}
%\caption{Three ways to map a $2$-simplex tightly to a tree: the image can
%be (1) a triod, (2) a path, or (3) a point. In (1) and (2) the
%image might have more vertices than are shown; also, each point inverse
%is a properly embedded arc, except in case (1) where one is a properly
%embedded triod.}
%\label{FigureSimplexMap}
%\end{figure}

The first step in building the 2-complex $X$ is to find the vertex set
$X^0$.  We will  build the vertex set using the $G$ action on the ends of
$T$ (note that even though $G$ only quasi-acts on $T$ it still honestly
acts on the ends).

If $G$ actually acts on $T$ then our construction gives for $X$ the
complex with $1$-skeleton the dual graph of $T$, with $2$-cells attached
around the vertices of $T$ so as to make the dual graph simply connected.
This picture should make the construction easier to follow.

\vfill 3.3. {\it Quasi\/{\rm -}\/edges}. \vfill

From the action at infinity we want to get some finite action. Each edge~$e$ of $T$ cuts $T$ into two ``sides'', which
are the two components of
$T-\interior(e)$, each a subtree of $T$. The end spaces of the two sides
of $e$ partition $\Ends(T)$ into an unordered pair of subsets denoted
$\E(e) = \{\C_1,\C_2\}$. Each of these $\C_1, \C_2$ is a {\it clopen} of
$\Ends(T)$ which means a subset that is both closed and open.

Generalizing this, we define a {\it quasi-edge} of $T$ to be a
decomposition of $\Ends(T)$ into a disjoint, unordered pair of clopens
$\E=\{\C_1,\C_2\}$, and $\C_1, \C_2$ are called the {\it sides} of the
quasi-edge $\E$. Although $G$ does not act on the set of edges
{\it a priori\/}, clearly $G$ acts on the set $\QE(T)$ of quasi-edges.

We now define the ``distortion'' of each quasi-edge $\E$ of $T$. This is a
positive integer $R(\E)$ which measures how far $\E$ is from
being a true edge. In particular, $R(\E)$ will equal $1$ if and only if
$\E = \E(e)$, as described above, for a unique edge~$e$.

Consider a clopen subset $\C$ of $\Ends(T)$. Let $\Hull(\C)$ denote the
convex hull of $\C$ in $T$, the smallest subtree of $T$ whose set of
accumulation points in $\Ends(T)$ is $\C$. Equivalently,
$\Hull(\C)$ is the union of all bi-infinite geodesics in $T$ whose
endpoints lie in $\C$. Under the inclusion $\Hull(\C)\inject T$ we may
identify $\Ends(\Hull(\C))$ with \pagebreak  $\C \subset \Ends(T)$.

Consider a quasi-edge $\E=\{\C_1,\C_2\}$ of $T$. Since each
vertex of $T$ has valence at least $3$, it follows that $V \subset
\Hull(\C_1) \union
\Hull(\C_2)$; for if there existed $v \in V - \bigl(\Hull(\C_1) \union
\Hull(\C_2)\bigr)$ then by convexity at least one of the three or more
components of
$T-v$ would be disjoint from both $\Hull(\C_1)$ and $\Hull(\C_2)$, and all
the
ends of $T$ reachable by that component would be disjoint from both $\C_1$
and
$\C_2$, contradicting that $\E=\{\C_1,\C_2\}$ is a quasi-edge. A similar
argument, together with connectivity of
$T$, shows that either $\Hull(\C_1) \intersect \Hull(\C_2) \ne \emptyset$
or the shortest path in $T$ connecting a point of $\Hull(\C_1)$
to a point of $\Hull(\C_2)$ is an edge. It
follows that there is at least one edge $e$ of $T$ such that $\bdy e
\intersect \Hull(\C_i)\ne \emptyset$ for each $i=1,2$; moreover, if
$\Hull(\C_1) \intersect \Hull(\C_2) \ne
\emptyset$ then there
are at least three such edges. Let $\Core(\E) = \Core(\C_1,\C_2)$ denote
the union of all such edges $e$; equivalently,
$$\Core(\E) = N_1(\Hull(\C_1)) \intersect N_1(\Hull(\C_2))
$$
where $N_1$ denotes the neighborhood of radius $1$ in $T$. Since
$N_1(\Hull(\C_i))$ is convex it follows that $\Core(\E)$ is a subtree
of $T$. In fact
$\Core(\E)$ is a finite subtree: if it were an infinite subtree then it
would accumulate on some end of~$T$, producing a point contained in the
closures of both $\C_1$ and $\C_2$, contradicting that $\E=\{\C_1,\C_2\}$
is a quasi-edge. We have also seen that
$\Core(\E)$ contains at least one edge---it is not a single point.
The number
$$R(\E) = \diam(\Core(\E))
$$
is therefore a positive integer, called the {\it quasi-edge distortion} of
$\E$. Given a positive integer $R$, if $R(\E) \le R$ then we
say that $\E$ is an {\it $R$-quasi-edge}.

It is  obvious that if $e$ is an edge, then the associated quasi-edge
$\E(e)$ satisfies $R=1$. Conversely, given a quasi-edge
$\E=\{\C_1,\C_2\}$ with $R(\E)=1$, it follows that $\Hull(\C_1)
\intersect \Hull(\C_2)$ is a single edge $e$. Moreover the two sides
of $e$ must be $\Hull(\C_1)$ and $\Hull(\C_2)$, and so $\E=\E(e)$.

As we saw above, a quasi-isometry takes quasi-edges to quasi-edges, and
now we investigate how the quasi-edge distortion is affected by a
quasi-isometry:

\proclaimtitle{Behavior of quasi-edge distortion under a quasi-isometry}
\specialnumber{14}\proclaim{Lemma} \label{LemmaQEDistortion}
For each $K,C$ there exists a constant $A$ such that if $\phi$ is a $K,C$
quasi\/{\rm -}\/isometry of $V=\Vertices(T)$ and if $\E=\{\C_1,\C_2\}$ is a
quasi\/{\rm -}\/edge then
$$d_\Haus\biggl(\phi\bigl(\Core(\C_1,\C_2)
\intersect V\bigr),\Core\bigl(\phi(\C_1),\phi(\C_2)\bigr)\biggr)
\le A 
$$
It follows that if $\E$ is an $R$\/{\rm -}\/quasi\/{\rm -}\/edge then
$\phi(\E)$ is a $KR+C+2A$ quasi\/{\rm -}\/edge 
\endproclaim

\demo{Proof}
There is a constant $A'$ depending only on $K,C$ such that if $\gamma$ is
a bi-infinite geodesic in $T$ with boundary $\bdy\gamma \subset \Ends(T)$
then $\phi(\gamma \intersect V)$ is within Hausdorff distance $A'$ of the
bi-infinite geodesic connecting the two points $\phi(\bdy\gamma)$. It
follows that
$$d_\Haus\bigl(\phi(\Hull(\C_i) \intersect V),\Hull(\phi(\C_i)\bigr) \le
A'
$$
and so
$$d_\Haus\bigl(\phi(N_1(\Hull(\C_i)) \intersect
V,N_1(\Hull(\phi(\C_i)))\bigr)
\le A'+K+C=A
$$
from which it follows that
\vglue12pt
\hfill $d_\Haus\bigl(\phi(\Core(\C_1,\C_2)
\intersect V),\Core(\phi(\C_1),\phi(\C_2))\bigr)
\le A 
$\enddemo
\vglue8pt
3.4. {\it Construction of the $2$\/{\rm -}\/complex $X$}.

\vglue6pt {\it The $0$\/{\rm -}\/skeleton of $X$.}
Consider the action of $G$ on $\QE(T)$. The $0$-skeleton $X^0$ consists
of the union of $G$-orbits of $1$-quasi-edges of $T$. By Lemma
\ref{LemmaQEDistortion} each element of $X^0$ is an $R$-quasi-edge where
$R=K+C+2A$, although perhaps not all $R$-quasi-edges are in $X^0$.
Clearly $G$ acts on $X^0$, because $X^0$ is a union of $G$-orbits of the
action of $G$ on the set of all quasi-edges $\QE(T)$. Define a map $f
\from X^0 \to \Vertices(T)$ by taking a quasi-edge $\E$ to any vertex in
$\Core(\E)$; by Lemma \ref{LemmaQEDistortion} this map is coarsely
$G$-equivariant. 

Since each edge $e \subset T$ determines a $1$-quasi-edge
$\E(e)$, since $\E(e) \in X^0$, and since $e=\Core(\E(e))$, it follows
that at least one of the two vertices of $e$ is in $f(X^0)$.

We claim that the cardinality $\abs{f^\inv(v)}$ is uniformly bounded
independent of $v$. It suffices to verify that there are boundedly many $R$-quasi-edges $\E$ such that $v \in
\Core(\E)$. To verify this, first note that the number of vertices of $T$ within distance $R$ of $v$ is at most
$1 + k(k-1)^{R-1}$, where $k$ is the maximum valence of a vertex of $T$.
It follows that there are boundedly many subtrees of $T$ of diameter
$\le R$ containing $v$. For each such subtree $\C$, there are boundedly
many components of $T-\C$. A quasi-edge with core $\C$ is determined by a
partition of the component set of $T-\C$ into two subsets, and there are
boundedly many such partitions. This proves the claim.

\vglue12pt{\it The $1$-skeleton of $X$.}
We now extend the $G$-set $X^0$ to a $1$-dimensional $G$-complex $X^1$,
by attaching edges to $X^0$ in two stages: first to make $X^1$
quasi-isometric to $T$, and second to extend the $G$ action.

In the first stage, two vertices $v,w \in X^0$ are attached by an edge if
$d(f(v),f(w)) = 0$ or $1$ in $T$, or if $d(f(v),f(w))=2$ in $T$ and the
vertex of $T$ between $f(v)$ and $f(w)$ is not in $f(X^0)$.  As noted
above,  each edge of $T$ contains at least one endpoint in $f(X^0)$.
Already at this stage the $1$-complex is connected, because for  any
$v,w \in X^0$, if $f(v)=f(w)$ then there is an edge from
$v$ to $w$, and  if $f(v) \ne f(w)$ then the path in $T$ from $f(v)$ to
$f(w)$ has at least every other vertex in $\image(f)$. Any further
attachment of cells of  dimension $\ge 1$ will preserve connectedness; we
use this fact without  comment \pagebreak from now on.

In the second stage, attach additional edges in a $G$-equivariant manner:
given vertices $v,w\in X^0$, if there exists $h \in G$ such that $h \cdot
v,h \cdot w$ are attached by a first stage edge, and if $v,w$ are not
already attached by a first stage edge, then $v,w$ are to be
attached by a second stage edge. This defines a connected $1$-complex
$X^1$, and an action of $G$ on $X^1$.

The map $f \from X^0 \to T$ is extended over $X^1$ by taking each edge $e
\subset X^1$ to the shortest path in $T$ connecting $f(v)$ to
$f(w)$, where $\bdy e = \{v,w\}$.

Putting the usual geodesic metric on $X^1$ where each edge has length
$1$, we claim the map $f \from X^1 \to T$ is a $G$-quasiconjugacy. To
prove this, we've already shown that $f \restrict X^0$ is coarsely
equivariant, and so it suffices to show that $f \restrict X^0$ is a
quasi-isometry. To see why, note first that each first stage edge of
$X^1$ maps to a path of length $\le 2$ in $T$, and each second stage edge
maps to a path of length $\le 2K+C$, where $K,C$ are quasi-isometry
constants for the $G$-quasi-action on $T$; this shows that $f$ is coarsely
lipschitz, in fact $f(d(v),d(w)) \le (2K+C) d(v,w)$ for $v,w \in
\Vertices(T)$.
To get the other direction, consider $v,v' \in X^0$. If $f(v)=f(v')$ then
$d(v,v')\le 1$. If $f(v) \ne f(v')$, let $f(v)=w_0 \to w_1 \to \cdots \to
w_k = f(v')$ be the geodesic in $T$ from $f(v)$ to $f(v')$; since at least
every other vertex $w_0,\ldots,w_k$ is in $f(X^0)$ it follows there is
an edge path in $X^1$ from $v$ to $v'$, consisting entirely of first
stage edges, of length $\le k$. In either case we've shown that $d(v,v')
\le d(f(v),f(v'))+1$.

We note some additional facts about the $1$-complex $X^1$:
\vglue3pt
\noindent\hskip16pt\hangindent=35pt\hangafter=1 (1)\  For each $w \in \Vertices(T)$ the set
$f^\inv(w) \intersect X^0$ has uniformly bounded cardinality.
\vglue3pt
\noindent\hskip16pt\hangindent=28pt\hangafter=1 (2)\    The vertices of $X^1$ have bounded valence.
\vglue3pt
\noindent\hskip16pt\hangindent=34pt\hangafter=1 (3)\   For each edge $e$ of $T$ and each $x \in
\interior(e)$, the set
$f^\inv(x)$ is a finite set of uniformly bounded cardinality
and diameter in $X^1$.
 \vglue3pt

Fact (1) was demonstrated earlier. Facts
(2) and (3) both follow from
Fact (1), together with the fact that for each edge $v
\stackrel{e}{\to} v'$ in $X^1$ the distance $d(f(v),f(v'))$ is uniformly
bounded, and the fact that an edge in $X^1$ is determined by its
endpoints.

\vglue6pt {\it The $2$-skeleton of $X$.}
We claim that there is a constant $B$ such that attaching 2-cells
along  all simple loops in $X^1$ of length $\le B$ results in a simply
connected 2-complex. First note that there is a $B$ such that all
isometrically embedded loops have length $\le B$; this follows
immediately from the fact that $X^1$ is quasi-isometric to a tree
(indeed it holds for any Gromov hyperbolic graph). Any loop is freely
homotopic to a concatenation of simple loops, and any simple loop is
freely homotopic to a concatenation of isometrically embedded loops.
Thus, once 2-cells are attached along all simple loops of length $\le
B$, all loops are freely null homotopic, which proves the resulting
2-complex, $X$, is simply connected.

The action of $G$ on $X^1$ clearly permutes the set of simple loops of
length $\le B$, and therefore extends to an action of $G$ on $X$.

It will be convenient, in what follows, to alter the cell-structure on
$X$ to obtain a simplicial complex. Since $X^1$ is already a simplicial
complex, this can be done by taking each 2-cell $\sigma$, introducing
a new vertex in the interior of $\sigma$, and connecting this new vertex
to each original vertex of $\sigma$, thereby cutting $\sigma$ into $b$
2-simplices where $b$ is the number of edges of $\bdy\sigma$. We put
a $G$-equivariant geodesic metric on $X$ so that each 2-simplex is
isometric to an equilateral Euclidean triangle of side length $1$.

Now we extend the map $f$ in a $G$-equivarian manner to obtain a map $f
\from X \to T$. This map is already defined on the $1$-skeleton of the
original cell-structure on $X$. For each original 2-cell $\sigma$ of
$f$, let $v(\sigma)$ be the new vertex in the interior of $\sigma$. Map
$v(\sigma)$ to any vertex in $f(\bdy\sigma)$, and map the new edges to
the unique geodesic in $T$ connecting the images of the endpoints.

%It is now easy to extend $f$ over
%the 2-simplices of $X$ so that the resulting map is tight.

The inclusion $X^1 \inject X$ is a quasi-isometry, and so the map $f \from
X \to T$ is a quasi-isometry. Clearly $f$ quasiconjugates the $G$ action
on $X$ to the original quasi-action on $T$.

%
%
%\paragraph{New end of proof}\quad
%
%\begin{itemize}
%\item First get rid of tightness; its not needed. In fact, simplify as
%much about the $2$-complex as possible.
%\item For cobounded case, all that we need is that it is
%locally finite, its quasi-isometric to a tree, and its isometry group is
%cobounded.
%\item For cobounded case, quote Dunwoody {\it directly}.
%\begin{itemize}
%\item First define tracks
%\item Then state Dunwoody's theorem to the $2$-complex with respect to its
%isometry group. (first we will have to equivariantly triangulate the
%$2$-complex, by starring each $2$-cell).
%\item Then apply it to get a quotient tree to which the isometry group
%descends.
%\item Check that the quotient map is a quasi-isometry.
%\end{itemize}
%\end{itemize}

\vglue 12pt 3.5. {\it Tracks}.
\vglue6pt
Since $G$ quasi-acts coboundedly on $T$, and since coboundedness is a
quasiconjugacy invariant, it follows that $G$ acts coboundedly on $X$.
Using this fact, the proof of Theorem~1 will be
finished quickly once we recall Dunwoody's tracks.

A {\it track} in a simplicial 2-complex $Y$ is a
$1$-dimensional complex $t$ embedded in $Y$ such that for each
2-simplex $\sigma$ of
$Y$, $t \intersect \sigma$ is a disjoint union of finitely many arcs,
each of which connects points in the interiors of two distinct edges of
$\sigma$. For each edge $e$ of $Y$ and each $x\in t\intersect e$, each
2-simplex $\sigma$ incident to $e$ therefore contains a component of
$t \intersect \sigma$ incident to $x$. A track $t \subset Y$ has a
{\it normal bundle} $p \from N(t) \to t$, consisting of a regular
neighborhood $N(t) \subset Y$ of $t$ and a fiber bundle $p\from N(t)\to
t$ with interval fiber, such that $p$ collapses each fiber to the unique
point where that fiber intersects $t$. If $Y$ is simply connected then
the complement of each track in $Y$ has two components; in particular,
the track locally separates, i.e.\ its normal bundle is orientable. A
track is {\it essential} if it separates $Y$ into two unbounded
components. Two tracks are {\it parallel} if they are ambient isotopic,
via an isotopy of $Y$ which preserves the skeleta.

In what follows, we shall assume that all tracks are {\it finite}.

\proclaimtitle{Tracks theorem \cite{Dun85}}
\specialnumber{15}\proclaim{Theorem} 
\label{TheoremTracks}
If $Y$ is a locally finite{\rm ,} simply connected{\rm ,} simplicial $2$-\/{\rm }\/complex
with cobounded isometry group{\rm ,} then there exists a disjoint union of
essential tracks $\tau = \disjunion_i \tau_i$ in $Y$ which is invariant
under the action of $\Isom(Y)$ such that the closure of each component of
$Y-\tau$ has at most one end.
 \endproclaim

Now we prove Theorem~1.

Apply Dunwoody's theorem to $X$ obtaining a disjoint union of tracks
$\tau =\disjunion_i\tau_i$. Consider the closure $A$ of a component of
$Y-\tau$. By Dunwoody's theorem, the set $A$ has at most one end. We claim
that in fact $A$ is bounded; this follows from the fact that $X$ is
quasi-isometric to a tree, by a standard argument which we now recall.

Suppose that $A$ is unbounded. Let $\Stab(A)$ be
the subgroup of $\Isom(X)$ that stabilizes $A$. Since $\Isom(X)$ acts
coboundedly on $X$ it follows that $\Stab(A)$ acts coboundedly on $A$.
Choose a sequence of points $x_0, x_1, x_2, \ldots \in A$, all in the same
orbit of $\Isom(A)$, such that $(x_i)$ escapes to infinity in $X$, and so
by passing to a subsequence we may assume that $(x_i)$ converges to some
end $\eta$ of $X$. Since $A$ is connected, its image under the
quasi-isometry $X \to T$ contains a ray converging to $\eta$, and so we
may assume that $x_0, x_1, x_2, \ldots$ lie on a quasigeodesic ray in $X$
converging to $\eta$. Choose $g_i \in \Stab(A)$ such that $g_i(x_i)=x_0$.
Since $X$ is quasi-isometric to a tree, and since $x_0,x_1,x_2,\ldots$
lie on a quasigeodesic ray converging to $\eta$, there exists
$R>0$ such that for all sufficiently large $i$ the compact set $\overline
N_R(x_i)$ separates $x_0$ from $\eta$. It follows that $C_i =
\overline N_R(x_i) \intersect A$ is a compact subset of $A$ separating
$x_0$ from $\eta$. Note that $C_0 = g_i(C_i)$ for all $i$. Let $U_i$ be
the component of $A-C_i$ containing $x_0$ and let $V_i$ be the component
limiting on $\eta$, so $U_i \ne V_i$. It follows that $\diam(U_i) \to
\infinity$ as $i \to \infinity$, and $\diam(V_i)=\infinity$ for all $i$.
Since $C_0$ has only finitely many complementary components, we may pass
to a subsequence so that $g_i(U_i)$ and $g_i(V_i)$ are constant, equal to
$U,V$ respectively. But then $U,V$ are distinct components of $A-C_0$,
each of unbounded  diameter. This shows that $A$ has $\ge 2$ ends, a
contradiction.

We now construct a map from $X$ to a tree $T'$, equivariant with respect
to $G$, in fact equivariant with respect to the entire isometry group of
$X$. Choose an equivariant
system of regular neighborhoods $N_i = N(\tau_i)$, and choose
equivariantly a fibration of $N_i$ by tracks parallel to $\tau_i$. The
tree $T'$ is the quotient of $X$ obtained by collapsing the closure of
each component of $X - \mathbold{\union} N_i$ to a point producing a vertex of $T'$,
and collapsing each of the parallel tracks in $N_i$ to a point producing
an edge of $T'$. The quotient map $X \to T'$ is equivariant, and it is
a quasi-isometry because the point inverse images are bounded.

This finishes the proof of Theorem~1.

\section{Application: Quasi-isometric rigidity for graphs\\ of coarse
$\PD(n)$ groups}
\label{SectionCoarsePDn}
4.1.  {\it Bass\/{\rm -}\/Serre theory}.
\vglue6pt

We review briefly graphs of groups, their Bass-Serre trees, and
associated topological spaces \cite{Ser80}, \cite{SW79}.

A {\it graph of groups} is a graph or $1$-complex $\Gamma$, together
with the following data: a vertex group $\Gamma_v$ for each vertex $v \in
\Vertices \Gamma$; an edge group $\Gamma_e$ for each $e \in
\Edges(\Gamma)$; and for each end $\eta$ of each edge $e$, with $v(\eta)$
the vertex incident to $\eta$, an injective edge-to-vertex homomorphism
$\gamma_\eta \from \Gamma_e \to \Gamma_{v(\eta)}$. The {\it fundamental
group} $\pi_1\Gamma$ can be defined topologically by first constructing a
{\it graph of spaces} associated to $\Gamma$, as follows. For each $v
\in\Vertices(\Gamma)$ choose a based $K(\Gamma_v,1)$ space $Y_v$; for each
$e \in \Edges(\Gamma)$ choose a based $K(\Gamma_e,1)$ space $Y_e$; and
for each end $\eta$ of each edge $e$ choose a base point preserving map
$f_\eta \from Y_e \to Y_{v(\eta)}$ inducing the homomorphism
$\gamma_\eta$. Let $Y$ be the quotient space
$$Y = \left(\coprod_{v \in \Vertices(\Gamma)} Y_v\right)
{\textstyle\coprod}
\left(\coprod_{e \in \Edges(\Gamma)} Y_e
\cross e\right)
\Biggm/ (x,v(\eta)) \sim f_{\eta}(x)
$$
where the indicated gluing is carried out for each edge $e$, each $x \in
Y_e$, and each end $\eta$ of $e$. The homotopy type of $Y$ is completely
determined independent of the choices of $Y_v$, $Y_e$, and $f_\eta$, and
the fundamental group of $Y$ is defined to be the fundamental group of
$\Gamma$.

The {\it Bass-Serre tree} of $\Gamma$ can also be defined topologically,
as follows. Define the {\it fibers} of $Y$ to be the images under the
above quotient of the vertex spaces $Y_v$ and the spaces $Y_e \cross x$,
$x \in \interior(e)$. In the universal cover $\wt Y$, the connected lifts
of the fibers of $Y$ are defined to be the fibers of $\wt Y$. The quotient
space of $\wt Y$ obtained by collapsing each fiber to a point is a tree
$T$ on which $\pi_1\Gamma$ acts, with quotient $\Gamma$. The graph of
groups structure on $\Gamma$ can be recovered using the vertex and edge
stabilizers of $T$ and the inclusion maps from edge stabilizers to vertex
stabilizers.

Let $\Gamma$ be a graph of groups. The universal cover $\wt Y$ equipped
with its quotient map $\wt Y \to T$ is an example of a ``tree of spaces''
for $\Gamma$, a concept which we now generalize. A {\it tree of spaces}
for $\Gamma$ consists of a cell complex $X$ on which $\pi_1\Gamma$ acts
properly by cellular automorphisms, together with a\break
$\pi_1\Gamma$-equivariant cellular map
$\pi \from X \to T$, such that the following properties hold:
\begin{itemize}
\item For each vertex $v$ of $T$, the set $X_v=\pi^\inv(v)$ is a connected
subcomplex of $X$ called the {\it vertex space} of $v$, and the
stabilizer group of $v$, $\Stab(v)=\{g\in\pi_1\Gamma\suchthat g \cdot v =
v\}$, acts properly on $X_v$.
\item For each edge $e$ of $T$ there is a connected cell complex $X_e$ on
which $\Stab(e)$ acts properly, called the {\it edge
space} of $e$, and there is a $\Stab(e)$-equivariant cellular map $i_e
\from X_e \cross e \to X$, such that $\pi \composed i_e$ equals the
projection map $X_e\cross e \to e$, and such that $i_e \restrict X_e
\cross \interior(e)$ is a homeomorphism onto $\pi^\inv(\interior(e))$
taking open cells to open cells.
\end{itemize}
We will sometimes identify $X_e$ with
$\pi^\inv(\midpoint(e))$. For each vertex $v$ of each edge $e$ of $T$,
the composition $X_e \approx X_e \cross
v \inject X_e \cross e \stackrel{i_e}{\longrightarrow} X$ has image contained in $X_v$
and therefore defines a cellular map $\xi_{ev}
\from X_e \to X_v$ called an {\it edge-to-vertex map}. Regarding $X_e =
\pi^\inv(\midpoint(e))$, the map $\xi_{ev}$ moves each point of $X_e$ a
bounded distance in $X$ to a point of $X_{v}$.

If $\Gamma$ is a finite graph of finitely generated groups, one can
construct a tree of spaces $\pi \from X \to T$ on which the action of
$\pi_1\Gamma$ is cobounded, by choosing the vertex and edge spaces so
that the action of the corresponding stabilizer is cobounded, for example
by taking Cayley graphs. Then one chooses cellular edge-to-vertex maps
$\xi_{ev} \from X_e \to X_{v}$, so that $\pi_1\Gamma$ acts equivariantly
on this data.  Define $X$ by gluing up the edge and vertex
spaces, that is:
$$X = \left(\coprod_v X_v\right) {\textstyle\coprod} \left(\coprod_e X_e
\cross e\right)
\Biggm/ (x,v) \sim \xi_{ev}(x) 
$$
The projection map $\left(\coprod_v X_v\right) {\textstyle\coprod}
\left(\coprod_e X_e\cross e\right) \to T$, which takes $X_v$ to $v$ and
$X_e \cross e$ to $e$ by projection, agrees with the gluings and
therefore defines the map $\pi \from X \to T$.

\vglue 12pt 4.2. {\it Geometrically homogeneous graphs of groups}.
\vglue6pt
A finite graph of finitely generated groups $\Gamma$ is
{\it geometrically homogeneous} if any of the following equivalent
conditions hold:
\vglue5pt
$\bullet$ $T$ has bounded valence;
\vglue5pt
$\bullet$ each edge-to-vertex injection $\gamma_\eta$ of $\Gamma$ has finite
index image;
\vglue5pt
$\bullet$ each edge-to-vertex injection $\gamma_\eta$ of $\Gamma$ is a
quasi-isometry;
\vglue5pt
$\bullet$ each edge-to-vertex map $\xi_{ev}$ of $X$ is a quasi-isometry;
\vglue5pt
$\bullet$ any two edge or vertex spaces in $X$ have finite Hausdorff
distance in~$X$.
\vglue5pt
\noindent
In the last three statements, we use any finitely generated word metric
on the edge and vertex groups, geodesic metrics on the edge and vertex
spaces, and a geodesic metric on $X$, on which the appropriate groups act
isometrically. Note that the first two statements are equivalent for any
finite graph of groups, regardless of whether the edge and vertex groups
are finitely generated. As proved in \cite{BK90}, if these
properties hold then the Bass-Serre tree $T$ satisfies a trichotomy:
either $T$ is bounded; or $T$ is {\it line-like} meaning that it is
quasi-isometric to a line; or $T$ is bushy. In the latter case we will
also say that the graph of groups $\Gamma$ is bushy.

Geometric homogeneity implies that all edge and vertex spaces of $X$, and
all edge and vertex groups of $\Gamma$, are in the same quasi-isometry
class. The converse does not hold, however: for a counterexample, take a
group $G$ having a monomorphism $\phi\from G \to G$ with infinite index
image, for instance a finite rank free group, and consider the HNN
amalgamation \pagebreak     $G*_\phi$.

4.3. {\it Weak vertex rigidity}.
\vglue6pt
 Let $\Gamma$ be a geometrically
homogeneous graph of groups with bushy Bass-Serre tree $T$, and choose a
tree of spaces $\pi\from X \to T$ for $\Gamma$.

Let $H$ be a group and let $(h,x) \mapsto h\cdot x$ be a
$K,C$-quasi-action of $H$ on $X$. We say that the quasi-action satisfies
{\it weak vertex rigidity} if there exists $R \ge 0$ such that for each
$h \in H$ and each vertex $v \in\Vertices(T)$ there is a vertex $v' \in
\Vertices(T)$ such that
$$h \cdot X_v \ceq{[R]} X_{v'} 
$$
Choosing one such $v'$ for each $v$ we obtain an induced $(K,C+2R)$
quasi-isometry $A_h\from T \to T$.

We claim that $h \mapsto A_h$ is a quasi-action of $H$ on $T$. To verify
this we must estimate the sup distance between $A_{hh'}$ and $A_h
\composed A_{h'}$. If we set
$$v' = A_{hh'}(v)
\qquad\hbox{then}\qquad
X_{v'} \ceq{[R]} hh' \cdot X_v \ceq{[C]} h \cdot (h' \cdot X_v) 
$$
And if we set
$$v_1 = A_{h'}(v), \quad v_2 = A_h(v_1)
$$
then
\begin{eqnarray*}
&&\phantom{h \cdot (h' \cdot X_v)}h' \cdot X_v \ceq{[R]} X_{v_1} \\
&&h \cdot (h' \cdot X_v) \dline{[KR+C]}  h \cdot X_{v_1}
\ceq{[R]} X_{v_2} 
\end{eqnarray*}
and so 
$$ 
X_{v_2} \ldline{[(R)+(KR+C)+(C)+(R)]}  X_{v'}
 S $$
showing that the sup distance between $A_{hh'}$ and $A_h \composed
A_{h'}$ is at most $KR+2C+2R$. 

In other words, any weakly vertex rigid quasi-action of a group $H$ on
the tree of spaces $X$ induces a quasi-action of $H$ on the
Bass-Serre tree~$T$. In this situation we can apply
Theorem~1, obtaining a quasiconjugacy $f \from
T' \to T$ from an isometric action of $H$ on a bounded valence, bushy tree
$T'$ to the induced quasi-action of $H$ on $T$.

When the original quasi-action of $H$ on $X$ is cobounded and proper,
evidently the induced quasi-action of $H$ on $T$ is cobounded, and so the
isometric action of $H$ on $T'$ is cobounded. Thus the quotient
$\Gamma'=T'/H$ may be regarded as a geometrically homogeneous graph of
groups with fundamental group $H$. In this situation it easily follows
that for each vertex $v'$ of $T'$, the stabilizer subgroup $\Stab_H(v')$
is finitely generated, and so we may construct a cobounded tree of spaces
$\pi'\from X' \to T'$ for the graph of groups~$\Gamma'$.

We now have most of the pieces in place for the following \pagebreak result:

\specialnumber{16}\proclaim{Proposition}
\label{PropGeomHomogQIRigidity}
Let $\Gamma$ be a geometrically homogeneous graph of groups{\rm ,}
with tree of spaces $X \stackrel{\pi}{\to} T$. Let $H$ be a
finitely generated group and let $(h,x) \mapsto h \cdot x$ be a weakly
vertex rigid{\rm ,} proper{\rm ,} cobounded quasi\/{\rm -}\/action of $H$ on $X$. Let $(h,v)
\mapsto h \cdot v$ be an induced quasi-action on $T$. Then there exists a
geometrically homogeneous graph of groups $\Gamma'$ with fundamental
group $H$ and with cobounded tree of spaces $X' \xrightarrow{\pi'}
T'${\rm ,} and there exists a coarsely commutative diagram
$$
\begin{array}{ccccc}
X' &\nhs\stackrel{F}{\lrar}\nhs&  X&\nhs\stackrel{\bar F}{\lrar} \nhs& X'
\\ 
\scrs{\pi'} &\nhs\nhs&\scrs{\pi} &\nhs\nhs&\scrs{\pi'} \\
T' &\nhs\stackrel{f}{\lrar}     \nhs           & T &\nhs\stackrel{\bar f}{\lrar}     \nhs        & T'
\end{array}
$$
in which all horizontal arrows are quasiconjugacies{\rm ,} $F$ and $\bar F$ are
coarse inverses{\rm ,} and $f$ and $\bar f$ are coarse inverses. As a
consequence{\rm ,} all of the following objects are quasi\/{\rm -}\/isometric\/{\rm :} the vertex
and edge groups of $\Gamma'${\rm ;} the vertex and edge spaces of
$X'${\rm ;} the vertex and edge spaces of $X${\rm ;} the vertex and edge groups of
$\Gamma$.
\endproclaim

{\it Proof}. Construct the quasiconjugacy $F$ as follows: choose
arbitrarily a point $x' \in X'$ and its image $x=F(x') \in X$. For
each point $x''\in X$ in the\break $H$-orbit of $x'$, choose arbitrarily $h \in
H$ so that $h \cdot x' = x''$, and define $F(x'')=h\cdot x$. For any other
point $x''' \in X$, choose arbitrarily a closest point $x''$ in the
$H$-orbit of $x'$, and define $F(x''')=F(x'')$. It is straightforward to
check that $F$ is an $H$-quasiconjugacy, and that the two maps $\pi
\composed F , f \composed \pi' \from X' \to T$ have bounded distance in
the sup norm. The coarse inverse $\bar F$ is similarly constructed.

To complete the proof it suffices to show that vertex spaces of $X$ and
$X'$ are quasi-isometric, and for this it suffices to show that for each
vertex $v' \in T'$ we have $F(X'_{v'}) \ceq{c} X_{fv'}$. The coarse
inclusion $F(X'_{v'}) \csubset{c} X_{fv'}$ follows from coarse
commutativity of the left-hand square, for $\pi'(F(X'_{v'}))$ is
contained in a uniformly bounded neighborhood of $fv'$, which implies in
turn that $F(X'_{v'}) \csubset{c} X_{fv'}$. Together with coarse
commutativity of the right-hand square we get $\bar F(X_{fv'})
\csubset{c} X'_{\bar f(fv')} \ceq{c} X'_{v'}$. Putting these together we
see that $F(X'_{v'}) \ceq{c} F(\bar F(X_{fv'}) \ceq{c} X_{fv'}$.
\hfill\qed\vglue10pt

4.4.  {\it Coarse  Poincar{\rm \'{\it e}}  duality spaces and groups}. \vglue4pt
Proposition~\ref{PropGeomHomogQIRigidity} has a strong
hypothesis, namely weak vertex rigidity. In this section we prove the
proposition that any bushy graph of coarse $\PD(n)$ groups with fixed $n$
satisfies weak vertex rigidity. Combined with
Proposition~\ref{PropGeomHomogQIRigidity} this will finish the proof of
Theorem~\ref{TheoremCoarsePDnGraphs}.

In \cite{FM98} it was proved that any bushy graph of $\Z$'s
satisfies weak vertex rigidity, using the coarse algebraic topology first
developed in \cite{FS96}. In \cite{FM00} the same
argument was generalized to certain bushy graphs of aspherical
$n$-manifold groups for fixed~$n$, with extra hypotheses added to ensure
that at no place in the argument did one leave the category of aspherical
manifolds, in order that coarse algebraic topology could still apply (these
extra hypotheses would be unnecessary if the Borel Conjecture were true).
With the theory of coarse \Poincare\ duality spaces developed by Kapovich
and Kleiner \cite{KK99}, the same argument can now be
further generalized. We begin by recalling some basic concepts from the
cohomology of groups \cite{Bro82}.

A group $G$ is {\it of type} FP if $\Z$ admits a finitely generated
projective resolution of finite length over $\Z G$; equivalently, there
is a finite CW complex which dominates any $K(G,1)$. The
{\it cohomological dimension} $\cd(G)$ of an FP group $G$ is the smallest
$n$ for which there exists a finitely generated projective resolution $0
\to P_n \to \cdots \to P_0 \to \Z \to 0$ over $\Z G$.

A group $G$ is {\it of type} FL if $\Z$ admits a finitely generated free
resolution of finite length over $\Z G$. A group $G$ has
{\it finite type} if there is a finite $K(G,1)$ CW complex. The
Eilenberg-Ganea Theorem (\cite{Bro82}, Theorem 7.1) shows that
$G$ has finite type if and only if $G$ is finitely presented and of type
FL; if this is the case, moreover, then there exists a finite $K(G,1)$ of
dimension $\max\{3,\cd(G)\}$.

An FP group $G$ is of type $\PD(n)$, $n \ge 1$, if the canonical
cohomology $H^i(G;\Z G)$ is $\Z$ when $i=n$ and is trivial otherwise; in
this case, $n=\cd(G)$. Mike Davis has examples of $\PD(n)$ groups in every
dimension $n \ge 4$ which are not of type FL and hence not   finitely presented
\cite{Dav98}.

Next we review the concepts of coarse \Poincare\ duality spaces from
\cite{KK99}, extending these concepts from
simplicial complexes to CW complexes.

Given a finite dimensional CW complex $X$, we define what it means for
$X$ to have {\it bounded geometry}. This means that  each point of $X$
touches a uniformly bounded number of closed cells, and loosely speaking
there is a uniformly finite set of models for the attaching maps of
cells. The formal definition uses induction on dimension, as follows. A
$1$-dimensional CW complex $X^1$ has bounded geometry if the valence of
$0$-cells is uniformly bounded; note that for each $R$ there are only
finitely many cellular isomorphism classes of connected subcomplexes of
$X^1$ containing $\le R$ cells. Suppose $X^{n+1}$ is an
$n+1$ dimensional CW complex whose $n$-skeleton $X^n$ has bounded
geometry, and assume inductively that for each $R$ there are only
finitely many cellular isomorphism classes of connected subcomplexes of
$X^n$ containing $\le R$ cells. Then $X^n$ has bounded geometry if there
exists an integer
$A > 0$ with the following properties:
 \vfil
\noindent\hskip16pt\hangindent=31pt\hangafter=1 $\bullet$ \enspace Each point of $X^n$ touches at most $A$ closed
cells of dimension
$n+1$.
\vfil
\noindent\hskip16pt\hangindent=31pt\hangafter=1 $\bullet$ \enspace  For each $n+1$ cell $e$ with attaching map
$\alpha_e \from S^{n} \to X^n$, the set $\image(\alpha_e)$ is a subcomplex of $X^n$ with
at most $A$ cells.
\vfil
\noindent\hskip16pt\hangindent=31pt\hangafter=1 $\bullet$ \enspace   Up to postcomposition by cellular
isomorphism, there are at most~$A$ different attaching maps $S^n \xrightarrow{\alpha_e} \image(\alpha_e)$, as
$e$ varies over all $n+1$~cells.

\pagebreak\noindent
It follows that for each $R$ there are only finitely many cellular
isomorphism classes of connected subcomplexes of $X^{n+1}$ with $\le R$
cells, completing the induction.

Given two bounded geometry CW complexes $X,Y$, suppose $f \from X \to Y$
is a {\it cellular map}, meaning that $f(X^n) \subset Y^n$ and the image
of each cell of $X$ is a subcomplex of $Y$. We say furthermore that
$f$ has bounded geometry if there is a uniformly finite set of models for
the maps from cells of $X$ to their images in $Y$; formalizing this along
the above lines is straightforward.

For example, the universal cover $\wt X$ of a finite CW complex $X$ has
bounded geometry. Also, if $f \from X \to Y$ is a cellular map between
finite CW complexes then any lift $\wt f \from \wt X \to \wt Y$ has
bounded geometry. $\phantom{\sum^2}$

Given a bounded geometry CW-complex $X$ and a subcomplex $L$, the
$i^{\rm th}$ combinatorial neighborhood $N_i(L)$ is a subcomplex of $X$ defined
inductively as follows: $N_0(L)=L$; and $N_{i+1}(L)$ is the union of all
closed cells intersecting $N_i(L)$. Intuitively, we imagine that $X$ has a
geodesic metric in which each $1$-cell is an arc of length $1$, each
cell has diameter $1$, and disjoint cells have distance $\ge 1$; with
such a geodesic metric, $N_i(L)$ would be the true metric neighborhood
about $L$ of radius $i$, and the inclusion $X^1 \subset X$ would be
a quasi-isometry. In fact one can pick a geodesic metric on $X$ having
only finitely many isometry types of connected subcomplexes of dimension
$\le R$ for each $R$, so that the inclusion $X^1 \subset X$ is a
quasi-isometry, and so if $B(L,r)$ denotes the metric neighborhood about
a subcomplex $L$ then we have $B(L,\frac{i}{K}-C) \subset N_i(L) \subset
B(L,Ki+C)$ for constants $K \ge 1, C \ge 0$. We are therefore free to
treat the neighborhood $N_i(L)$, for quasi-isometric purposes, as a metric
neighborhood.

Suppose $X$ is a bounded geometry CW complex. We say that $X$ is
{\it uniformly contractible} if for each $R$ there exists $R'=R'(R) > R$
such that each subcomplex $L$ with $\diam(L) \le R$ is contractible to a
point in $N_{R'}(L)$ (here and in the sequel, ``distance'' measurements
in $X$ such as $R$ and $R'$ are implicitly assumed to be nonnegative
integers).

We say that the CW chain complex $C_*(X)$ is {\it uniformly acyclic} if
for every $R$ there exists $R'=R'(R)>R$ such that for every subcomplex $L
\subset X$ with $\diam(L)\le R$ the inclusion $L\to N_{R'}(L)$ induces
the trivial map on reduced homology.

Let $C^*_c(X)$ be the compactly supported CW cochain complex, let
$n$ be the topological dimension of $X$, and suppose we are given an
{\it $n$-dimensional augmentation}, which means a surjective
homomorphism $C^n_c(X)\xrightarrow{\alpha}\Z$ with the property that the
augmented sequence
$$C^0_c(X) \to \cdots \to C^{n-1}_c(X) \to C^n_c(X) \xrightarrow{\alpha}
\Z
$$
is a cochain complex. We say that the augmented cochain complex
$(C^*_c(X),\alpha)$ is {\it uniformly acyclic} if there exists $R_0$, and
for each $R \ge R_0$ there exists $R'=R'(R) > R$, such that the following
hold: 
 
\noindent\hskip16pt\hangindent=31pt\hangafter=1 $\bullet$ \enspace  for all vertices $v$ and all $R \ge R_0$ the image of the
restriction map
$$H^*_c(X,\overline{X - N_R(v)}) \to H^*_c(X,\overline{X-N_{R'}(v)})
$$
maps isomorphically to $H^*_c(X)$ under the restriction map
$$H^*_c(X,\overline{X-N_{R'}(v)}) \to H^*_c(X) 
$$
 
\noindent\hskip16pt\hangindent=31pt\hangafter=1 $\bullet$ \enspace  the induced map $\alpha^* \from H^n_c(X) \to \Z$ is an isomorphism.
\vglue4pt
\noindent\hskip16pt\hangindent=31pt\hangafter=1 $\bullet$ \enspace  for all $i \ne n$, $H^i_c(X) \approx 0$.
\vglue4pt

Following \cite{KK99}, a bounded geometry, uniformly
acyclic CW complex $X$ is a {\it coarse $\PD(n)$ space} if there exist
chain maps
$$C_*(X) \xrightarrow{P} C^{n-*}_c(X) \xrightarrow{\overline P} C_*(X)
$$
such that the maps $\bar P \composed P$, $P \composed \bar P$ are each
chain homotopic to the identity via respective chain homotopies
$$C_*(X) \xrightarrow{\Phi} C_{*+1}(X),
\quad C^*_c(X)\xrightarrow{\bar\Phi} C^{*-1}_c(X)
$$
and there exists a constant $D\ge 0$ such that for all cells $\sigma$
of $X$, the supports of each of $P(\sigma)$, $\bar P(\sigma)$,
$\Phi(\sigma)$, $\bar\Phi(\sigma)$ lie in $N_D(\sigma)$.

\specialnumber{17}\proclaim{Lemma}
\label{LemmaAugmentation}
If $X$ is an $n$\/{\rm -}\/dimensional{\rm ,} bounded geometry {\rm CW} complex{\rm ,} then $X$ is a
coarse $\PD(n)$ space if and only if there exists an $n$\/{\rm -}\/dimensional
augmentation $\alpha \from C^n_c(X) \to \Z$ such that $(C^*_c(X),\alpha)$
is uniformly acyclic.
\endproclaim

\demo{Proof}
The ``if'' direction comes from \cite{KK99}. For the
``only if'' direction, take $\alpha$ to be the composition of the
duality map $\overline P \from C^n_*(X) \to C_0(X)$ with the ordinary
augmentation $C_0(X) \to \Z$.
\enddemo

The main fact we will need about coarse $\PD(n)$ spaces is the coarse
separation theorem of \cite{KK99}, which says
that a uniformly proper embedding of a coarse $\PD(n-1)$ space
in a coarse $\PD(n)$ space coarsely separates the target into exactly two
deep components.

We define a {\it coarse $\PD(n)$ group} to be a group $G$ which
quasi-acts properly and coboundedly on some  
coarse $\PD(n)$ space $X$ of topological dimension $n$; equivalently, $G$
is quasi-isometric to $X$. As a shorthand, a  
coarse $\PD(n)$ space of topological dimension $n$ is called a
{\it good coarse $\PD(n)$ space}. Thus, a coarse $\PD(n)$ group is one
which is quasi-isometric to a good coarse $\PD(n)$ space.

For example, if $G$ is a $\PD(n)$ group of finite type then $G$ is a
coarse $\PD(n)$ group. If $n=2$ this follows from the Eckmann-M\"uller
theorem which implies that $G$ is the fundamental group of a closed,
aspherical surface (see \cite[p.\ 223]{Bro82}). If $n\ge 3$,
the Eilenberg-Ganea theorem (\cite[Th.\ 7.1]{Bro82}) says
that $G$ has an $n$-dimensional, finite $K(G,1)$ space $K$, and the action
of $G$ on the universal cover of $K$ demonstrates that $G$ is a coarse
$\PD(n)$ group.

The coarse $\PD(n)$ property of a group is by definition a quasi-isometry
invariant. In particular it follows that groups which are virtually
$\PD(n)$ of finite type are coarse $\PD(n)$.

\vglue 12pt 4.5. {\it Bushy graphs of coarse $\PD(n)$ groups}.
\vglue9pt

Here is a restatement of
Theorem~\ref{TheoremCoarsePDnGraphs}:

\vglue9pt {\elevensc Theorem 18.} 
{\it Given $n\ge 0${\rm ,} if $\Gamma$ is a finite graph of coarse $\PD(n)$ groups with bushy
Bass\/{\rm -}\/Serre tree{\rm ,} and if $G$ is a finitely generated group quasi\/{\rm -}\/isometric
to $\pi_1\Gamma${\rm ,} then $G$ is the fundamental group of a graph of groups
with bushy Bass\/{\rm -}\/Serre tree{\rm ,} and with vertex and edge groups
quasi\/{\rm -}\/isometric to those of $\Gamma$.}
\vglue4pt
For expository purposes we will describe the proof in detail under the
following:
\vglue12pt {\it Standing assumption}. Each vertex and edge group {\it ACTS},
properly and coboundedly, on a good coarse $\PD(n)$ space.
\vglue8pt \noindent
In other words, while the coarse $\PD(n)$ property for groups only
requires the group to quasi-act, our standing assumption replaces
this with the stronger requirement that the group acts. Once the proof is
complete under this assumption, we sketch the changes needed to
cover the general case.

To start with we construct a tree of spaces $\pi \from X \to T$
by choosing vertex spaces and edge spaces which are good coarse $\PD(n)$
spaces on which the respective vertex and edge groups act properly and
coboundedly by cellular automorphisms; also, choose bounded geometry
edge-to-vertex maps. This construction is possible, first of all, because
of the standing assumption above. Secondly, we need to make use of the
fact that in a geometrically homogeneous graph of groups, each
edge-to-vertex map is a quasi-isometry; and so in our present situation
where the edge and vertex spaces are bounded geometry    it follows that the edge-to-vertex maps can each be
moved a bounded distance in the sup norm to a bounded geometry map. This can be
done equivariantly, because of the standing assumption, and by gluing we
obtain the required $\pi\from X\to T$.

Note that $X$ is a bounded geometry  $n+1$
complex. This follows because the constants for the bounded geometry of vertex and edge spaces and the constants for
the bounded geometry of edge-to-vertex maps are uniformly bounded.

Given a bi-infinite line $L \subset T$, the subcomplex $X_L = \pi^\inv(L)$
is called the {\it hyperplane} in $X$ lying over $L$. Each hyperplane is
a bounded geometry $n+1$ complex, and is uniformly
properly embedded in~$X$. Lemma~\ref{LemmaCoarsePDn} will show 
that each $X_L$ is a good coarse $\PD(n+1)$ space.

\specialnumber{19}\proclaim{Proposition}
\label{PropWeakVertexRigidity}
Let $P$ be a good{\rm ,} coarse $\PD(n+1)$ space. For every $K \ge 1${\rm ,} $C \ge
0${\rm ,} and every proper function $\rho \from[0,\infinity)\to [0,\infinity)$
there exists $A \ge 0$ such that if $f \from P \to X$ is a
$(K,C,\rho)$\/{\rm -}\/uniformly proper embedding{\rm ,} then there is a hyperplane
$X_{L}$ such that $d_\Haus(f(P),X_{L}) \le A$.
\endproclaim

\phantom{snow}
\vglue-24pt
This is the exact analogue of Theorem~7.3 in \cite{FM00}, and
will be proved below.

\vglue4pt {\it Proof\/{\rm :} Theorem~{\rm 18} reduces to
Proposition~{\rm \ref{PropWeakVertexRigidity}}}. The proof follows
\cite[\S7.2, Step 2]{FM00} very closely; here is a sketch.

Combining with Proposition~\ref{PropGeomHomogQIRigidity}, it remains to
show that any quasi-action on $X$ satisfies weak vertex rigidity. For
this we need only show the following: if $X' \to T'$
is another tree of spaces associated to a geometrically homogeneous graph
of coarse $\PD(n)$ groups, then for each $K,C$ there exists $A$ such that
if $f\from X \to X'$ is a $K,C$-quasi-isometry and if $v\in T$ is a vertex
then there is a vertex $v' \in T'$ such that $d_\Haus(f(X_v),X'_{v'})\le
A$. First note that there exists a proper function $\rho \from
[0,\infinity)\to [0,\infinity)$, depending only on
$K,C$, such that if $X_L$ is any hyperplane in $X$ then $f \restrict X_L$
is a $\rho$-uniformly proper embedding of $X_L$ into $X'$. Now choose
three bi-infinite lines $L_1,L_2,L_3$ in $T$ whose intersection is $v$.
Applying Proposition~\ref{PropWeakVertexRigidity} there are hyperplanes
$X'_{L'_1}, X'_{L'_2}, X'_{L'_3}$ which are uniformly Hausdorff close to
$f(X_{L_1}), f(X_{L_2}), f(X_{L_3})$ respectively. The set $f(X_v)$ must
therefore be uniformly Hausdorff close to the coarse intersection of the
hyperplanes
$X'_{L'_1}, X'_{L'_2}, X'_{L'_3}$. The three lines $L'_1, L'_2, L'_3$ must
intersect pairwise in rays, and these three rays have infinite Hausdorff
distance in $T'$; it follows that $L'_1 \intersect L'_2 \intersect L'_3$
is a vertex $v'$ of $T'$. This implies that the coarse intersection of
$X'_{L'_1}, X'_{L'_2}, X'_{L'_3}$ is $X'_{v'}$, proving that $f(X_v)$ is
uniformly Hausdorff close to $X'_{v'}$. 
\hfill\qed\vglue6pt

We turn now to the proof of Proposition~\ref{PropWeakVertexRigidity}.

The reader who is interested in most real life examples can refer to
\cite[Th.\ 7.3]{FM00} where the proof is given under the
following special circumstances: the vertex and edge groups of $\Gamma$
are fundamental groups of closed, smooth, aspherical manifolds in a
category $\C$ which is closed under finite covers, and which has the
property that any homotopy equivalence between manifolds in
$\C$ is homotopic to a diffeomorphism. For example: euclidean manifolds;
hyperbolic manifolds; irreducible, nonpositively curved symmetric spaces;
solvmanifolds; and nilmanifolds.

The proof of Proposition~\ref{PropWeakVertexRigidity} follows the same
outline as Theorem~7.3 of \cite{FM00}, but with many changes of
details needed to accommodate coarse $\PD(n)$ groups.

Pick a topologically proper embedding of $T$ in an open disc $D$. For
each component $U$ of $D-T$ there is a line $L(U)$ in $T$ such that the
projection homeomorphism $L(U) \cross 0 \to L(U)$ extends to a
homeomorphism of pairs $(L(U) \cross [0,\infinity), L(U) \cross 0) \homeo
(\closure U,L(U))$. Extend the CW structure on $T$ to a CW structure on
$D$, by using the product structure on $L(U) \cross [0,\infinity) \homeo
\closure U$ for each $U$, where $L(U)$ has the CW structure it inherits
from $T$ and $[0,\infinity)$ has the CW structure where each interval
$[i,i+1]$ is a $1$-cell. Note that with this CW-structure, $D$ is a
coarse $\PD(2)$ space.

The tree of spaces $\pi \from X \to T$ is an example of what might be
called a ``coarse fibration'', and we now extend this to a ``coarse
fibration'' $\pi\from Y\to D$. The fibers of $Y$ will be good coarse
$\PD(n)$ spaces isomorphic to the fibers of~$X$.

For each cell $c$ of $D$ we define a fiber $Y_c$ as follows. If $c
\subset T$ then we simply take $Y_c=X_c$. If $c \subset U$ for some
component $U$ of $D-T$, it follows that $c = c_1 \cross c_2$ for cells
$c_1\subset L(U)$ and $c_2 \subset [0,\infinity)$, and we define $Y_c$ to
be a disjoint copy of $X_{c_1}$.

For any two cells $c \supset d$ of $D$ we define an attaching map
$\eta_{dc}\from Y_c \to Y_d$ as follows. If $d,c \subset T$ we simply take
$\eta_{dc} = \xi_{dc}$. Otherwise we have $d,c \subset\overline U$ for
some component $U$ of $D-T$. Let $d=d_1 \cross d_2, c=c_1 \cross c_2$
where $d_1,c_1$ are cells of $L(U)$ and $d_2,c_2$ are cells of
$[0,\infinity)$. If $d_1=c_1$ then $Y_d, Y_c$ are disjoint copies of
$Y_{d_1}$ and we take $\eta_{dc}$ to be a disjoint copy of the identity
map on $Y_{d_1}$. Otherwise we take $\eta_{dc}$ to be a disjoint copy of
the map~$\xi_{d_1c_1}$.

Note that for cells $c \supset d \supset e$ of $D$ we have $\eta_{ce} =
\eta_{de} \composed \eta_{cd}$. We therefore have a well-defined
quotient space
$Y$ as follows:
$$Y = \coprod_c (Y_c \cross c) \Biggm/ (x,p) \sim (\eta_{cd}(x),p)
$$
where the pair $c,d$ varies over all cells $c \supset d$ of $D$, the
point $x$ varies over $Y_c$, and $p$ varies over all points of
$\bdy c$ that lie on $d$. The disjoint union of the projection maps $Y_c
\cross c \to c$ gives a map $\coprod_c(Y_c \cross c) \to D$, and the
latter is consistent with all gluings, thereby defining the map $\pi
\from Y \to D$. We identify $Y_c$ with $\pi^\inv(x)$ for a chosen point
$x$ in the interior of $c$, and we call this the {\it fiber} over $c$.
Once we have a metric in place it is evident that any two fibers in $Y$
have finite Hausdorff distance.

Here are a few facts about $Y$. First, since the base $D$ and fibers
$Y_c$ are bounded geometry cell complexes with uniform bounds, and since
the gluing maps $\eta_{cd}$ have bounded geometry with uniform bounds, it
follows that $Y$ has bounded geometry.  Also, the inclusion map from any fiber
$Y_c$ into $Y$ is uniformly proper, the inclusion map $X \inject Y$ is
uniformly proper, and the inclusion map from any hyperplane $X_L$
into $Y$ is uniformly proper.

\specialnumber{20}\proclaim{Lemma}
\label{LemmaCoarsePDn}
$Y$ is a good coarse $\PD(n+2)$ space. Each hyperplane $X_L$ is a good
coarse $\PD(n+1)$ space with uniform bounds independent of the line~$L$.
\endproclaim

Accepting this lemma for the moment, we now have:

\demo{{P}roof\/{\rm :} Proposition {\rm \ref{PropWeakVertexRigidity}} reduces to
Lemma {\rm \ref{LemmaCoarsePDn}}} The proof follows very closely
\cite[\S7.2, Step 1]{FM00}, which itself follows very closely
\cite[Prop.\ 4.1, Steps 1 and 2]{FM98}. Here is a sketch.

Applying the coarse Jordan separation theorem of
\cite{KK99} together with Lemma~\ref{LemmaCoarsePDn}
it follows that there exists $R > 0$ independent of $L$ such that $Y -
N_R(f(P))$ has exactly two components $A,A'$ which are {\it deep},
meaning that each of $A,A'$ contains points arbitrarily far from $f(P)$.

Since $f(P) \subset X$ and since each component of $Y-X$ is deep, it
follows that each component of $Y-X$ is coarsely contained in one of
$A,A'$, and each of $A,A'$ coarsely contains at least one component of
$Y-X$. We may therefore find two components $U,U'$ of $Y-X$ coarsely
contained in $A,A'$ respectively, such that $U,U'$ are adjacent in $Y-X$,
meaning that $\overline U \intersect \overline U' = \pi^\inv(E_0)$
for some edge $E_0$ of $T$. Letting $E_n$ denote the neighborhood of
radius $n$ about $E_0$ in $T$, a simple inductive argument shows that
$E_n$ contains an embedded edge path $\gamma_n$ of length $2n+1$ centered
on $E_0$ such that $\pi^\inv(\gamma_n)$ is coarsely contained in $f(P)$
with uniform coarse containment constant; if this were not so then as in
\cite[Prop.\ 4.1, Step 1]{FM98}   one can find a path in
$Y-N_R(f(P))$ connecting
$A$ to $A'$, a contradiction. It follows that $f(P)$ coarsely contains
some hyperplane $X_{L}$. A packing argument as in \cite[Prop.\ 4.1, Step 2]{FM98}
 shows that $X_{L}$ is coarsely contained in
$f(P)$.\hfill\qed
\enddemo

\demo{Proof  of Lemma {\rm \ref{LemmaCoarsePDn}}}
The proofs for $Y$ and for $X_L$ are exactly the same. We give the proof
for $Y$, using the $E_0$ spectral sequence in compactly supported CW
cohomology for the coarse fibration $\pi \from Y \to D$.

Here is some notation. For each cell $d$ of a CW complex $Z$ we let $d^*$
denote the corresponding basis element of the compactly supported CW
cochain complex $C^*_c(Z)$. The coboundary operator
$\cbdy \from C^i_c(D) \to C^{i+1}_c(D)$ can be written as
$$\cbdy b^* = \sum_c n_{cb} c^*, \quad\hbox{for $b$ an $i$-cell of $D$}
$$
where $c$ ranges over cells of dimension $i+1$ and, letting $m_c \from
S^{i} \to X^{i}$ be the attaching map for the cell $c$, the integer
$n_{cb}$ is the degree of the map $(S^{i},\emptyset) \xrightarrow{m_c}
(X^i, X^i - \interior(b))$.

Each cell in the CW complex $Y$ has the form $b \cross e$ for some cell
$b$ of $D$ and some cell $e$ of $Y_b$, and the corresponding basis
element of $C^*_c(Y)$ can therefore be written $b^* \cross e^*$. To
understand the coboundary operator of $Y$ we look at the $E_0$ term of
the spectral sequence.

The filtration of $D$ by its skeleta $D^0 \subset D^1 \subset D^2=D$
determines a filtration $Y^0 \subset Y^1 \subset Y^2 = Y$ with
$Y^i=\pi^\inv(D^i)$, which in turn determines an $E_0$ spectral sequence
for $H^*_c(Y)$ as follows. Over the CW complex $D$ we have a coefficient
bundle $C^*_c(\Y)$: the coefficients over a cell $b$ of $D$ are
$C^*_c(Y_b)$; and for cells $b \subset c$ the pullback  homomorphism is
$\eta_{cb}^\ast \from C^*_c(Y_b) \to C^*_c(Y_c)$. We then have
$$E_0^{ij} = C^i_c(D,C^j_c(\Y)) 
$$
The coboundary operator of each coefficient complex $C^*_c(Y_b)$
determines a map
$$d \from E_0^{ij} \to E_0^{i,j+1}
$$
given by $d(b^* \cross e^*) = b^* \cross \cbdy e^*$. The coboundary
map of
$C^*_c(D)$ together with the restriction maps of the coefficient bundle
determine a map
$$\cbdy \from E_0^{ij} \to E_0^{i+1,j}
$$
given by
$$\cbdy(b^* \cross e^*) = \sum_{b\subset c}   n_{cb} (c^*\times \eta^*_{cb}(e^*) )
$$
We have an isomorphism
$$C^k_c(Y) \homeo \bigoplus_{i+j=k} C^i_c(D,C^j_c(\Y))
$$
and under this isomorphism the cochain map for $C^*_c(Y)$ corresponds to
the map 
$$d + (-1)^k \cbdy \from \bigoplus_{i+j=k} C^i_c(D,C^j_c(\Y)) \to
\bigoplus_{i+j=k+1} C^i_c(D,C^j_c(\Y)) 
$$
The $E_1$ term of the spectral sequence is given by
$$E_1^{ij} = C^i_c(D,H^j_c(\Y)) =
\left\{ \begin{array}{ll}
\Z & \hbox{if $0 \le i \le 2$, $j=n$} \\
0  & \hbox{otherwise}
\end{array}\right.
$$
and the spectral sequence collapses at the $E_2$ term with
$$E_2^{ij} =
\left\{ \begin{array}{ll}
\Z & \hbox{if $i=2$, $j=n$} \\
0 &  \hbox{otherwise}.
\end{array}\right.
$$
This shows at least that $Y$ has the correct cohomology for a coarse
$\PD(n+2)$ space: $H^{n+2}_c(Y)=\Z$ and $H^i_c(Y)=0$ for $i \ne n+2$. But
we must still verify the uniformity criteria, and for this we look more
closely at the representation of $C^*_c(Y)$ using the double complex
$E_0$.

The CW complex $Y$ has topological dimension $n+2$, and so to verify that
$Y$ is a good coarse $\PD(n+2)$ space using Lemma~\ref{LemmaAugmentation}
we must construct an $n+2$ dimensional augmentation $\alpha$ for
$C^*_c(Y)$ and prove that $(C^*_c(Y),\alpha)$ is uniformly acyclic. Since
$D$ is a good coarse $\PD(2)$ space, there is a \hbox{2-dimensional}
augmentation
$\alpha_D$ making $C^*_c(D)$ uniformly acyclic; in fact, we simply define
$\alpha_D(\beta)=\sum_b \<\beta,b\>$, summed over 2-cells $b$ of $D$.
For each cell $d$ of $D$, since $Y_d$\break is a good coarse \pagebreak $\PD(n)$ space
then by Lemma~\ref{LemmaAugmentation} there is an $n$-dimensional
augmentation $\alpha_d$ making $C^*_c(Y_d)$ uniformly acyclic (with
constants independent of $d$). Define $\alpha \from C^{n+2}_c(Y) \to \Z$
as follows. A basis element of $C^{n+2}_c(Y)$ has the form $b^* \cross
e^*$, for $b$ a $2$-cell of $D$ and $e$ an $n$-cell of $Y_b$. We define
$$\alpha(b^* \cross e^*) = \alpha_D(b^*) \alpha_b(e^*) = \alpha_b(e^*) 
$$
It is clear from the spectral sequence that $\alpha$ determines an
isomorphism\break $\alpha^* \from H^{n+2}_c(Y) \to \Z$. Moreover, each class in
$H^n_c(Y_b) \approx \Z$ is represented by a cochain whose support is
contained in an arbitrary ball of uniformly bounded radius in $Y$, and so
the same is true for each class in $H^{n+2}_c(Y)$.
\vglue2pt
It remains to show, given an arbitrary coboundary $\gamma \in C^k_c(Y)$,
that $\gamma = \cbdy \rho$ for some $\rho \in C^{k-1}_c(Y)$ such that
$\supp(\rho)$ is contained in a uniformly bounded neighborhood of
$\supp(\gamma)$. 

\vglue12pt {\it Case} 1: $k<n+2$.  Using the double complex $E_0$ we may write
$$\gamma = \gamma^{0k} \oplus \gamma^{1,k-1} \oplus \gamma^{2,k-2} \in
E_0^{0k} \oplus E_0^{1,k-1} \oplus E_0^{2,k-2} 
$$
Since $\gamma$ is a coboundary, there exists $\rho = \rho^{0,k-1} \oplus
\rho^{1,k-2} \in E_0^{0,k-1} \oplus E_0^{1,k-2}$ with
\begin{eqnarray*}
d \rho^{0,k-1} &= &\gamma^{0k} \\
(-1)^{k-1}\,\delta \rho^{0,k-1} + d \rho^{1,k-2} &=& \gamma^{1,k-1} \\
(-1)^{k-1}\,\delta \rho^{1,k-2} &=& \gamma^{2,k-2} 
\end{eqnarray*}
Each $(C^*_c(Y_b),\alpha_b)$ is uniformly acyclic with uniform constants,
and so in each stalk of the coefficient bundle $C^*_c(\Y)$ we can choose
$\rho^{0,k-1}$ to have support contained in a uniformly bounded
neighborhood of $\supp(\gamma^{0k})$. Subtracting $(d + (-1)^{k-1} \cbdy)
\rho^{0,k-1}$ from $\gamma$ it suffices to assume that $\rho^{0,k-1}=0$
and so we have $d\rho^{1,k-2} = \gamma^{1,k-1}$. Repeating the argument,
we can choose $\rho^{1,k-2}$ to have support contained in a uniformly
bounded neighborhood of $\supp(\gamma^{1,k-1})$.

\vglue12pt {\it Case} 2: $k=n+2$.  We have $\gamma=\gamma^{2,n}$, and since
$\gamma$ is a coboundary it follows that $\alpha(\gamma)=0$. We may write
$\gamma = \sum_b \gamma_b$ summed over 2-cells $b$ of $D$ where
$\gamma_b \in C^n_c(Y_b)$, and so
$$0 = \alpha(\gamma) = \sum_b \alpha_b(\gamma_b) 
$$
We may reduce to the case that $\gcf\{\alpha_b(\gamma_b)\} = 1$. To see
why, for each $b$ choose $\gamma'_b \in C^n_c(Y_b)$ with
$\alpha_b(\gamma'_b)=1$ so that $\supp(\gamma'_b)$ is in a uniform
neighborhood of $\supp(\gamma_b)$. Then replace $\gamma_b$ with
$\alpha_b(\gamma_b) \cdot\gamma'_b$ by subtracting off the coboundary of
something in $C^{n-1}_c(Y_b)$ whose support is contained in a uniform
neighborhood of $\supp(\gamma_b)$. After this replacement, all
coefficients of the cocycle $\gamma$ are divisible by
$\gcf\{\alpha_b(\gamma_b)\}$, and carrying out this division produces a
cocycle with integer coefficients, finishing the reduction.

Define $\gamma_D \in C^2_c(D)$ to be the projection of $\gamma$ to $D$,
namely $\gamma_D = \sum_b \alpha_b(\gamma_b) b^*\break\in C^2_c(D)$ and so
$\alpha_D(\gamma_d)=0$. Since $(C^*_c(D),\alpha_D)$ is uniformly acyclic
we can find $\rho_D \in C^1_c(D)$ with $\cbdy\rho_D=\gamma_D$ such that
$\supp(\rho_D) \subset N_R(\supp(\gamma_D))$ for some uniform $R$.

We shall now construct, uniformly for each $1$-cell $a \in
\supp(\rho_D)$, an element $\rho_a \in C^{n}_c(Y_a)$ such that
$\alpha_a(\rho_a) = \<\rho_D,a\>$. To do this, for each inclusion $a
\subset b$ of a $1$-cell into a 2-cell in $D$ choose a cochain
map $\bar\eta_{ab} \from C^*_c(Y_a) \to C^*_c(Y_b)$ so
that $\eta_{ba} \composed \bar\eta_{ab}$ and $\bar\eta_{ab} \composed
\eta_{ba}$ are uniformly chain homotopic to the respective identities.
For each $1$-cell $a \in \supp(\rho_D)$ and
2-cell $b \in \supp(\gamma_D)$, choose an alternating sequence of
2-cells and $1$-cells
$$b=b_0, a_0, b_1, a_1, \ldots, b_r, a_r = a
$$
so that $r$ is uniformly bounded ($r \le xR+y$, with $x,y$ depending only
on $D$), and so that
$a_i\subset b_i$, $i=0,\ldots,r$ and $a_{i-1} \subset b_i$,
$i=1,\ldots,r$. Define
$$\gamma_a^b = \eta_{b_ra_r} \composed \bar\eta_{a_{r-1}b_r} \composed
\cdots \composed \bar\eta_{a_0b_1} \composed \eta_{b_0a_0} (\gamma_b)
\in C^n_c(Y_a) 
$$
Noting that $\eta_{ba}$ and $\bar\eta_{ab}$ each commute with
augmentations, it follows that
$$\alpha_a(\gamma^b_a) = \alpha_b(\gamma_b)
$$
and so for each $a$ we have $\gcf_b\{\alpha_a(\gamma^b_a)\} = 1$. This
implies that we can find a linear combination in $C^n_c(Y_a)$ of
the set of cocycles
$\{\gamma^b_a
\suchthat b \subset \supp(\gamma_D)\}$ so that the resulting cocycle
maps to $1$ under $\alpha_a$. Multiplying this cocycle by $\<\rho_D,a\>$
we obtain the desired $\rho_a \in C^n_c(Y_a)$. Note that $\supp(a^*
\cross \rho_a)$ is contained in a uniform neighborhood of
$\supp(\gamma)$. 

Consider now the $(n+2)$-cocycle $\hat\gamma = \gamma - \cbdy(\sum_a a^*
\cross\rho_a)$, and note that $\alpha_b(\hat\gamma_b)=0$ for each
2-cell $b$ of $D$. The cocycle $\hat\gamma_b \in C^n_c(Y_b)$ is
therefore the coboundary of some $\rho_b \in C^{n-1}_c(Y_b)$ whose
support is contained in a uniform neighborhood of $\hat\gamma_b$. It
follows that
$$\gamma = \cbdy\left(\sum_{a \subset \supp(\rho_D)} a^*
\cross
\rho_a + \sum_{b \subset \supp(\gamma)} b^* \cross \rho_b \right)
$$
and the $n+1$-cocycle inside the parentheses has support contained in a
uniform neighborhood of $\supp(\gamma)$.

This completes the proof of Lemma~\ref{LemmaCoarsePDn},
Proposition~\ref{PropWeakVertexRigidity}, and
Theorem~18, under the standing assumption.

\vglue12pt{\it Removing the standing assumption}. For the general proof we
choose, for each vertex or edge group, a vertex and edge space on which
the corresponding group quasi-acts properly and coboundedly, and which is
a good coarse $\PD(n)$ space. We can still glue the vertex and edge
spaces together to  \pagebreak form a tree of spaces $\pi \from X \to T$, and all the
key geometric properties hold: $X$ is bounded geometry, and the hyperplanes $X_L$ are still good coarse $\PD(n+1)$
spaces. Moreover, we can extend this to $\pi \from Y \to D$ as before,
where $Y$ is a good coarse $\PD(n+2)$ space.

What we have lost is that the group $\pi_1\Gamma$ no longer acts on $X$.
However, the key point is that $\pi_1\Gamma$ still quasi-acts on
$X$, properly and coboundedly. It follows that $\pi_1\Gamma$, as well as
any group quasi-isometric to $\pi_1\Gamma$, is quasi-isometric to $X$.
Moreover, the proof that $\pi \from X \to T$ satisfies weak vertex
rigidity is exactly the same, and so
Theorem~18 is proved.
\enddemo
 \vglue5pt
{\it Remark on {\rm ``}\/coarse fibrations\/{\rm ''}}. With the proper notion of a
coarse fibration, the proof above   generalizes to show that any
``coarse fibration'' over a good coarse $\PD(m)$ base space with good
coarse $\PD(n)$ fiber is a good coarse $\PD(m+n)$ space. This   implies
a generalization of a theorem of Bieri \cite{Bie72} and Johnson-Wall
\cite{JW72} saying that any extension of a $\PD(m)$
group by a $\PD(n)$ group is $\PD(m+n)$; the generalized statement
  replaces ``$\PD(n)$ group'' by ``coarse $\PD(n)$ group''.

\section{Application: Actions on Cantor sets}
\label{SectionCantorActions}
%\input{QTPartI_Cantor.tex}
% QTPartI_Cantor.tex

We shall prove Theorem \ref{TheoremQC} about uniformly quasiconformal
actions on Cantor sets, Theorem 7 about uniform
quasisimilarity actions on the $n$-adic rational numbers, and the
corollaries to these theorems.

\vglue 12pt 5.1. {\it Uniformly quasiconformal actions}.
\vglue6pt

Theorem \ref{TheoremQC} is an almost immediate corollary of
Theorem~1, once we review the results of
\cite{Pau96} which prove the equivalence of various notions
of quasiconformality for homeomorphisms between boundaries of Gromov
hyperbolic spaces, all of which are equivalent to being the extension of
a quasi-isometry; see \cite{Pau96} and
\cite{Pau95} for the history of these various notions.

First we review a notion of quasiconformality which was developed for
rank one symmetric spaces by Pansu \cite{Pan89a} and generalized
to word hyperbolic groups by Paulin \cite{Pau96}. Let $X,Y$ be
proper, geodesic hyperbolic metric spaces with cobounded isometry groups.
Let $\H(\bdy X,\bdy Y)$ denote the space of homeomorphisms from $\bdy X$
to $\bdy Y$ with the compact open topology. For each compact set $\K
\subset \H(\bdy X,\bdy Y)$, a homeomorphism $f \from \bdy X \to \bdy Y$ is
said to be {\it $\K$-quasiconformal} if for each isometry $\alpha \from X
\to X$ there exists an isometry $\beta \from Y \to Y$ such that
$\bdy\beta\composed f \composed \bdy\alpha \in \K$.

Next we review a notion of quasiconformality based on cross-ratio.
Consider a proper, geodesic, hyperbolic metric space
$X$, and let $\bdy^4 X$ denote the space of ordered $4$-tuples of pairwise
distinct points in
$\bdy X$. The
{\it cross-ratio} on $\bdy X$ is the function which assigns to each
$(a,b,c,d)\in\bdy^4 X$ the positive number
\begin{eqnarray*}
&&\hskip-24pt[a,b,c,  d]  \\
 &&= \exp\left( \frac{1}{2}
\sup_{{x_i \to a, \, y_i \to b \atop  z_i \to c, t_i \to d}}
     \liminf_{i \to \infinity} \bigl(
           d(x_i,t_i) - d(t_i,z_i) + d(z_i,y_i) - d(y_i,x_i)
     \bigr)
\right)
\end{eqnarray*}  
where the sequences $x_i, y_i, z_i, t_i$ are in $X$ and the convergence to
$a,b,c,d$, respectively, takes place in the Gromov compactification
$\overline X = X \union \bdy X$. When $X=T$ is a tree, no matter how
these sequences are chosen the expression in the parentheses eventually
takes on the same constant value, which is described simply as follows.
Consider the six geodesic lines in $T$ determined by taking the points
$a,b,c,d$ in pairs: $\overline{ab}$, $\overline{ac\vphantom{d}}$, etc.
Then we have
$$\abs{\log[a,b,c,d]} = \max\left\{ d(\overline{ab},\overline{cd}),
d(\overline{ac\vphantom{d}},\overline{bd}), d(\overline{ad},\overline{bc})
\right\} 
$$

Given two proper, geodesic, Gromov hyperbolic metric spaces $X,Y$ with
cobounded isometry group, and given a proper,
increasing function $I \from [0,\infinity) \to [0,\infinity)$, a
homeomorphism $f \from \bdy X \to\bdy Y$ is {\it $I$-quasimobius} if
for all $(a,b,c,d)\in\bdy^4 X$ we have
\begin{eqnarray*}
[a,b,c,d] &\le& I([fa,fb,fc,fd]) \\
\left[fa,fb,fc,fd\right] &\le& I([a,b,c,d]) 
\end{eqnarray*}
As noted in \cite{Pau96}, it turns out that one can restrict
attention to functions $I$ of the form $I(r) = a (\sup\{1,r\})^\kappa$.

\proclaimtitle{\cite{Pau96}}
\specialnumber{21}\proclaim{Theorem} 
\label{TheoremPaulin}
Let $X,Y$ be proper{\rm ,} geodesic{\rm ,} Gromov hyperbolic metric spaces with
cobounded isometry group. Given $f \from \bdy X \to \bdy Y${\rm ,} the
following are equivalent\/{\rm :}
\begin{itemize}
\ritem{(1)} There exists a compact $\K \subset \H(\bdy X,\bdy Y)$ such that $f$
is $\K$\/{\rm -}\/quasiconformal.
\ritem{(2)} There exists $K \ge 1, C \ge 0$ such that $f$ is the continuous
extension of some $K,C$ quasi\/{\rm -}\/isometry $X \mapsto Y$.
\ritem{(3)} There exists some proper{\rm ,} increasing $I \from [0,\infinity) \to
[0,\infinity)$ such that $f$ is $I$\/{\rm -}\/quasimobius.
\end{itemize}
Moreover{\rm ,} these equivalences are uniform\/{\rm :} for example{\rm , (1)} implies {\rm (2)}
uniformly in the sense that for any compact $\K \subset \H(\bdy X,\bdy Y)$
there exists\break $K \ge 1, C \ge 0$ such that if $f \from \bdy X \to \bdy Y$
is $\K$\/{\rm -}\/quasiconformal then $f$ is the continuous extension of a $K,C$
quasi\/{\rm -}\/isometry\/{\rm ;} similarly for the other five implications.
\endproclaim

Define a group action $A \from G \cross \bdy X \to \bdy X$ to be
{\it uniformly quasiconformal} if the collection of homeomorphisms $A_g
\from \xi \mapsto A(g,\xi)$ satisfies any of the criteria of
Theorem~\ref{TheoremPaulin} uniformly; for instance, there exists a
proper, increasing $I \from [0,\infinity)\to[0,\infinity)$ such that each
$A_g$ is $I$-quasimobius.

\demo{Proof  of Theorem {\rm \ref{TheoremQC}}}
Let $T$ be a bushy tree of bounded valence, and let $A \from G \cross
B \to B$ be a uniformly quasiconformal action of $G$ on $B=\bdy
T$ whose induced action on the triple space of $B$ is cocompact. For each
$g \in G$ apply Theorem \ref{TheoremPaulin} to extend the map
$A_g \from B \to B$ to a $K,C$ quasi-isometry $\alpha_g\from T \to T$,
with $K,C$ independent of $g$. For any $g, g' \in G$, since $\alpha_g
\composed\alpha_{g'}$ is a $K^2, KC + C$ quasi-isometry, the following
standard lemma shows that $\alpha_g\composed \alpha_{g'}$ is a bounded
distance in the sup norm from $\alpha_{gg'}$, proving that $g\mapsto
\alpha_g$ is a quasi-action of $G$ on $T$:

\specialnumber{22}\proclaim{Lemma} 
\label{LemmaQuasiComposition}
For any proper{\rm ,} Gromov hyperbolic metric space $X$ and any $K \ge 1, C
\ge 0$ there exists $A \ge 0$ such that if $f,g \from X \to X$ are
$K,C$ quasi\/{\rm -}\/isometries of $X$ whose boundary extensions $\bdy f, \bdy g
\from \bdy X \to \bdy X$ are identical{\rm ,} then the {\rm sup} distance between $f$
and $g$ is at most $A$. 
\endproclaim

Since the action of $G$ on the triple space of $B$ is cocompact, it
follows that the quasi-action $g \mapsto A_g$ is cobounded. Applying
Theorem 1 to the quasi-action $g
\mapsto A_g$ we finish the proof of Theorem \ref{TheoremQC}.
\enddemo

\demo{Proof  of Corollary {\rm \ref{CorollaryQC}}} Let $\alpha' \from G
\to \Isom(T')$ be the action produced by Theorem~\ref{TheoremQC}, so that
$B$ is identified quasiconformally with $\bdy T'$. By coboundedness of
the action, it follows that the quotient graph of groups $\Gamma=T'/G$ is
finite. Since $T'$ has bounded valence, $\Gamma$ has finite index
edge-to-vertex injections. Given a subgroup $H$ of $G$, suppose the
action of $H$ on the space of distinct pairs in $B$ has precompact
orbits. This implies that the action of $H$ on $T'$ has bounded orbits,
which implies that $H$ fixes some vertex of
$T'$.
\enddemo

{\it Remark}. In Section 2.2 of \cite{Pau96} there is
another
intrinsic notion of a quasiconformal structure on $\bdy X$, based on
moduli of annuli, which says roughly that a homeomorphism is
quasiconformal if it stretches the moduli of annuli by a uniform amount;
the modulus stretching function $\phi$ is then a measure of uniformity of
the homeomorphism. However, this measure of uniformity does not agree
uniformly with the ones used in Theorem~\ref{TheoremPaulin}; for example,
as noted at the end of \cite[\S2.2]{Pau96}, isometries of $X$
do not even extend with uniform quasiconformality modulus $\phi$. Perhaps
with some tinkering this notion of quasiconformality could be made to
agree uniformly with the others.

\demo{{R}emark} It could be interesting to have an intrinsic notion of
a conformal structure on a Cantor set $B$ so that a homeomophism
$B\xrightarrow{f} \bdy T$, where $T$ is a bounded valence, bushy tree,
determines a conformal structure on $B$ with the property that conformal
homeomorphisms are the same as boundary extensions of isometries of $T$.
Paulin's paper discusses conformal structures on boundaries of Gromov
hyperbolic metric spaces, but these do not behave well on Cantor set
boundaries or other boundaries which do not have a sufficiently rich set
of arcs; for example, using that notion of conformality one can find
conformal maps of $\bdy T$ whose extension to $T$ can be made isometric
outside of a finite subtree of $T$, but which cannot be made isometric on
all of $T$. Instead, it might be possible to use quasi-edges to come up
with a better behaved notion of conformality on $B$.
\enddemo

5.2. {\it Uniform quasisimilarity actions on $n$-adic Cantor sets}.
\vglue6pt

First we review some basic notions about maps of metric spaces. Recall
from the introduction that given a metric space $X$, a
$K$-quasisimilarity is a bijection $f \from X \to X$ such that
$$\frac{d(f\zeta,f\omega)}{d(\zeta,\omega)} \Biggm/
\frac{d(f\xi,f\eta)}{d(\xi,\eta)}
\le K, \quad\hbox{for all } \zeta \ne \omega, \xi \ne \eta \in X.
$$
A $K$-bilipschitz homeomorphism is a bijection $f\from X\to Y$ between
metric spaces such that
$$\frac{1}{K} d(\xi,\eta) \le d(f\xi,f\eta) \le K d(\xi,\eta), \quad
\xi,\eta \in X.
$$ 
If $K$ is unspecified in either of these definitions then $f$ is called
simply a {\it quasisimilarity} or {\it bilipschitz}, respectively.

Given a compact subinterval $[a,b] \subset (0,\infinity)$, a bijection
$f\from X \to X$ is {\it $[a,b]$-bilipschitz} if
$$\frac{d(f\xi,f\eta)}{d(\xi,\eta)} \in [a,b], \quad\hbox{for all } 
\xi\ne\eta\in X.
$$
We use the notation $r \cdot [a,b]$ for the interval $[ra,rb]$.

As observed in \cite{FM99}, each $[a,b]$-bilipschitz
homeomorphism of a metric space $X$ is a $K$-quasisimilarity with
$K=\frac{b}{a}$, and each $K$-quasisimilarity is $[a,b]$-bilipschitz for
some interval $[a,b]$ with $\frac{b}{a} \le K^4$.

\bigskip

Next we review results from \cite{FM98} concerning the
connection between $\Q_n$ and trees.

Let $T_n$ be the homogeneous directed tree with one edge pointing towards
each vertex and $n$ edges pointing away, and put a geodesic metric on
$T_n$ where each edge has length $1$. In the end space $\bdy T_n$, a
Cantor set, there is a unique end denoted $-\infinity$ which is the limit
point of any ray in $T_n$ obtained by starting at a vertex in $T_n$ and
travelling backwards against the direction of edges. There is a
{\it height map} $\height\from T_n\to \R$ which takes each directed edge
of $T_n$ isometrically to a directed interval $[i,i+1]$ with integer
endpoints; the map $\height$ is unique up to postcomposition by
translation of $\R$ by an integer amount, and we shall fix once and for
all a choice of $\height$. For each $i \in \Z$ the set $\height^\inv(i)$
is called the {\it level set} of height $i$ in $T_n$. A {\it coherent
line} in $T_n$ is a continuous section of the height function $T_n \to
\R$. There is a one-to-one correspondence between the set of coherent
lines and the set $\bdy T_n - \{-\infinity\}$, where the coherent line
$\ell$ corresponds to the point $\xi$ if $\bdy\ell =
\{-\infinity,\xi\}$; we denote $\ell=\ell_\xi$.

There is a homeomorphism between $\bdy T_n - \{-\infinity\}$ and $\Q_n$,
determined uniquely up to isometry of $\Q_n$ by the property that for
all $\xi \ne \eta \in \Q_n$, if the vertex at which $\ell_\xi$ and
$\ell_\eta$ diverge from each other has height $h$ then
$$d(\xi,\eta) = n^{-h} 
$$
We get an explicit picture of $T_n$ as follows.

A {\it clopen} in a topological space is a subset which is both closed
and open. The {\it separation} of a subset $U$ in a metric space $X$ is
defined to be
$$\sep(U) = \inf\{d(x,y) \suchthat x \in U, y \in X - U\} 
$$
Note that if $U$ is a clopen in $\Q_n$ then $\sep(U) > 0$.

A {\it clone} in $\Q_n$ is a certain kind of clopen, defined as follows.
Given an integer $h$, define an equivalence relation on $\Q_n$ where
$\xi,\eta \in \Q_n$ are equivalent if $\xi_i=\eta_i$ for all $i \le h$;
the equivalence classes are called {\it clones of height $h$}, and each
of them is a clopen in $\Q_n$. Thus, there is one clone of height $h$ for
every sequence $\omega=(\omega_i)_{i=-\infinity}^h$ in $\Z/n$
satisfying the property that $\omega_i$ is eventually zero as $i\to
-\infinity$; the integer $h$, which specifies that the domain of the
sequence $\omega$ is the interval $(-\infinity,h]$, is also called the
height of the sequence $\omega$.

There is a one-to-one, height preserving correspondence between vertices
of $T_n$ and clones of $\Q_n$, where the vertex $v$ corresponds to the
clone $U_\omega$ if and only if the positive endpoints of the coherent
lines passing through $v$ are exactly the points of the clone $U_\omega$;
we write $v=v_\omega$ and $\omega=\omega_v$ to emphasize this
correspondence. Note also that the structure of $T_n$ as an oriented tree
corresponds to the inclusion structure among clones: there is a directed
edge $v_\omega \to v_{\omega'}$ if and only if $U_{\omega'} \subset
U_\omega$ and no other clone is nested strictly between $U_{\omega}$ and
$U_\omega$; moreover this occurs if and only if $\omega$ has some height
$h$, $\omega'$ has height $h+1$, and $\omega_i=\omega'_i$ for all $i \le
h$. 

\medskip

Now we describe certain isometries and quasi-isometries of $T_n$
and their effect on $\Q_n$.

A {\it height translation} of $T_n$ is an isometry $f \from T_n \to T_n$
with the property that $h_0=\height(f(v))-\height(v)$ is constant for $v
\in T_n$; the constant $h_0$ is called the {\it height translation
length} of $f$, and
$f$ is called more specifically an {\it $h_0$-height translation}.
Height translations form a subgroup of the full isometry group of $T_n$,
in fact they are exactly the orientation preserving isometries of $T_n$.
Height translations of $T_n$ are related to similarities of
$\Q_n$ as follows:

\specialnumber{23}\proclaim{Proposition}
\label{PropHeightTranslations}
Continuous extension defines an isomorphism between the height translation
group of $T_n$ and the similarity group of $\Q_n$. The extension to $\Q_n$
of an $h_0$\/{\rm -}\/height translation of $T_n$ is an $n^{-h_0}$ similarity of
$\Q_n$. In particular{\rm ,} the height preserving isometry group of $T_n$
corresponds to the isometry group of $\Q_n$.  
\endproclaim

We also need a quasification of Proposition~\ref{PropHeightTranslations}.
A quasi-isometry $f \from T_n\to T_n$ is called a {\it coarse height
translation} if there exists $A\ge 0$ such that $\height(f(v))-\height(v)$
varies over some subinterval $[m,m+A] \subset \R$ of length $\le A$ as $v$
varies over $T_n$. Any value $h_0=\height(f(v))-\height(v) \in [m,m+A]$ is
called a
{\it coarse height translation length} of $f$. We shall incorporate the
coarseness constant $A$ and the height translation length $h_0$ into the
terminology by referring to $f$ as an {\it $A$-coarse $h_0$-height
translation}. 

The quasification of Proposition \ref{PropHeightTranslations} says roughly
speaking that continuous extension defines a ``uniform'' isomorphism
between the coarse height translation group of $T_n$ and the
quasisimilarity group of $\Q_n$:

\specialnumber{24}\proclaim{Proposition}
\label{PropCoarseHeightTranslations}  \hskip-8pt
For each $K' \ge 1$ there exist constants $K \ge 1, C
\ge 0${\rm ,} $A \ge 0${\rm ,} $R \ge 1$ with the following properties\/{\rm :}
\begin{itemize}
\item  Each $K'$\/{\rm -}\/quasisimilarity $F \from
\Q_n\to\Q_n$ extends to a $(K,C)$ quasi\/{\rm -}\/isometric $A$\/{\rm -}\/coarse height translation $f \from T_n
\to T_n${\rm ;}
\item  If $f$ is a coarse $h_0$\/{\rm -}\/height translation then $F$ is $n^{-h_0}
\cdot [R^\inv,R]$\/{\rm -}\/bilipschitz.
\end{itemize}
Conversely{\rm ,} for each $K \ge 1, C \ge 0, A \ge 0$ there exists $K'
\ge 1$ such that the continuous extension of each $(K,C)$\/{\rm -}\/quasi\/{\rm -}\/isometric
$A$\/{\rm -}\/coarse height translation of $T_n$ is a $K'$\/{\rm -}\/quasisimilarity of
$\Q_n$. 
\endproclaim

\demo{Proof  of Theorem~{\rm 7}} Fix a uniform
quasisimilarity action of a group $G$ on $\Q_n$, $n\ge 2$, and suppose
that the induced action of $G$ on the space of distinct doubles of $\Q_n$
is cocompact. 

We must produce a bilipschitz
homeomorphism $\Q_n \to \Q_p$ which conjugates the $G$ action on $\Q_n$ to
a similarity action on $\Q_p$.

\vglue6pt {\it Step} 1: {\it A topological conjugacy.}
Using Proposition~\ref{PropCoarseHeightTranslations} together with
Lemma~\ref{LemmaQuasiComposition}, the uniform quasi-similarity action of
$G$ on $\Q_n$ extends to a quasi-action of $G$ on $T_n$ by uniformly
coarse height translations. Moreover, coboundedness of the $G$-action on
the double space of $\Q_n$ translates to coboundedness of the
$G$-quasi-action on $T_n$.

Apply Theorem~1 to get a new tree $T'$, a
cobounded isometric action of $G$ on $T'$, and a quasiconjugacy $f \from
T'\to T_n$ with continuous extension $F \from \bdy T' \to \bdy T_n$.
Since the quasi-action of $G$ on $T_n$ fixes the end $-\infinity$, the
action of $G$ on $T'$ fixes the end $F^\inv(-\infinity)$ which we also
denote as $-\infinity$. Orient $T'$ away from $-\infinity$, and so the
$G$-action respects this orientation. The downward valence at each vertex
of $T'$ with respect to this orientation equals $1$.

Now we make some simplifications to $T'$.

First, we may assume:
\vglue8pt
  ($*$)  $T'$ has no valence~1 vertices.

\pagebreak\noindent 
For if there is a valence~1 vertex, incident to an edge $e$, we may
equivariantly collapse all edges in the orbit $G\cdot e$ producing a new
tree $T''$ on which $G$ acts, and a $G$-equivariant collapsing map $T'
\mapsto T''$. The collapsing map $T' \mapsto T''$ is a quasi-isometry
because each component of the forest $\mathbold{\union} G \cdot e$ is a graph of
uniformly finite size: in fact, each component is a star graph with
$k$ valence~1 vertices for some constant~$k$ and with one vertex of
valence~$k$. Replacing $T'$ with $T''$ reduces the number of vertex
orbits, of which there are finitely many. Continuing inductively,
eventually ($*$) is established. 

At this stage, the quotient graph $T'/G$ is a connected directed graph,
with inward valence~1 at each vertex, and outward valence at least~1. It
follows that $T'/G$ is a circle with one or more vertices, all of whose
directed edges agree with some global orientation on the circle.

For the second simplification, we may assume:
\vglue8pt
\noindent\hskip16pt\hangindent=31pt\hangafter=1 $\bullet$ \enspace $T'$ has exactly one
orbit of vertices and one orbit of edges, or in other words the
directed graph $T'/G$ is a circle with one vertex and one edge.
\vglue8pt\noindent
If this
is not so, consider any edge $e$ of $T'$, pointing toward a vertex
$\bdy_+ e$, and pointing away from a vertex $\bdy_- e$ in a different
orbit than $\bdy_+ e$. Since $\bdy_+ e \ne \bdy_- e$, and since no other
edge in the orbit $G \cdot e$ points toward $\bdy_+ e$, it follows that
the components of $\mathbold{\union} G \cdot e$ are graphs of uniformly finite size.
Collapsing all edges in the orbit $G \cdot e$ produces a tree $T''$ on
which $G$ acts and a $G$-equivariant, quasi-isometric collapse $T' \to
T''$. Replacing $T'$ by $T''$ and continuing inductively, eventually ($**$)
is established.

With ($*$) and ($**$) established, $T'$ is an oriented tree on which $G$
acts, transitively on vertices and on edges, with downward valence $1$ and
constant upward valence $p$ for some integer $p \ge 2$. In other words, we
may identify $T'=T_p$ as oriented trees. Extension of the $G$-action on
$T_p$ gives a similarity action on $\Q_p$. Since the quasiconjugacy $f
\from T_p \to T_n$ takes the $-\infinity \in \bdy T_p$ to
$-\infinity\in\bdy T_n$, it follows that $f$ induces a topological
(indeed quasiconformal) conjugacy
$F\from\Q_p\to\Q_n$ between  $G$-actions.

\vglue8pt
 {\it Step {\rm 2:} $F$ is bilipschitz}. We will reduce this statement to
Claim~\ref{ClaimDiamsSeps} which describes the effect of $F$ on clones of
$\Q_p$; in Step~3 we will prove  Claim~\ref{ClaimDiamsSeps}.
\vglue8pt
Let $\U_h$ be the set of height~$h$ clones in $\Q_p$, and let $F(\U_h) =
\{F(U) \suchthat U \in \U_h\}$ be the set of image clopens in $\Q_n$.

\specialnumber{25}\proclaim{Claim}
\label{ClaimDiamsSeps}
There exists $\lambda > 1$ and an interval $[a,b] \subset (0,\infinity)$
such that for all $h\in\Z$ and all $U \in \U_h${\rm ,}
$$\diam(F(U)), \, \sep(F(U)) \in \lambda^{-h} \cdot [a,b] 
$$
\endproclaim

Using this claim we show that $F$ is \pagebreak bilipschitz.

First we evaluate $\lambda$. Since each element of $\U_h$ subdivides into
$p$ disjoint elements of $\U_{h+1}$, it follows that each element $F(U)
\in F(\U_h)$ is partitioned into $p$  disjoint elements of $F(\U_{h+1})$;
moreover, by Claim~\ref{ClaimDiamsSeps} the diameters of these partition
elements shrink by a factor of $\lambda$ relative to the diameter of
$F(U)$, up to a bounded multiplicative error. A simple calculation shows
that the Hausdorff dimension of $\Q_n$ equals $\log_\lambda p$. But we
know that the Hausdorff dimension of
$\Q_n$ equals $1$, proving that
$$\lambda = p 
$$

Consider now two points $x,y \in \Q_p$, and let $U \in \U_h$ be the
smallest clone of $\Q_p$ containing $x,y$, so
$$d_{\Q_p}(x,y) = p^{-h} 
$$
We also have:
$$d_{\Q_n}(Fx,Fy) \le \diam(F(U)) \le b \lambda^{-h} = b p^{-h} 
$$
And, since $x,y$ are contained in {\it distinct} level $h+1$ clones
$U_x, U_y$ we have
$$d_{\Q_n}(Fx,Fy) \ge \sep(F(U_x)) \ge a \lambda^{-h} = a p^{-h}
$$
proving that
$$a d_{\Q_p}(x,y) \le d_{\Q_n}(Fx,Fy) \le b d_{\Q_p}(x,y)
$$
and so $F$ is bilipschitz.

\vglue12pt {\it Step {\rm 3:} Proof of Claim} \ref{ClaimDiamsSeps}. The similarity
action of $G$ on $\Q_p$ induces a stretch homomorphism $s \from G \to
(0,\infinity)$ defined by the property that $g$ is an $s(g)$-similarity
for each $g \in G$. Let $G_0 = \kernel(s)$, and so $G_0$ is the subgroup
of $G$ acting isometrically on $\Q_p$. The quotient group $G/G_0 \approx
\image(s)$ is infinite cyclic. Choose an infinite cyclic
subgroup $Z=\<z\>$ of $G$ such that the projection from $Z$ to
$G/G_0$ is an isomorphism. To prove the claim we shall show:
\begin{itemize}
\item[(1)]
The action of $G_0$ on $\Q_n$ is uniformly
bilipschitz: there exists $R \ge 1$ such that each $g \in G_0$ is a
$R$-bilipschitz homeomorphism of $\Q_n$.
\item[(2)]
There exists $\lambda>1$ and an interval $[a_0,b_0] \in
(0,\infinity)$ such that for each integer $k$, the action of $z^k$ on
$\Q_n$ is a $\lambda^{-k}\cdot [a_0,b_0]$-bilipschitz homeomorphism.
\end{itemize}
To see why this proves the claim, note first that $G$ acts transitively on
vertices of $T_p$ and so $G_0$ acts transitively on each level of $T_p$,
in particular on level zero. This shows that $G_0$ acts
transitively on the clones $\U_0$ of $\Q_p$. Since $F \from \Q_p \to
\Q_n$ is $G_0$-equivariant it follows that $G_0$ acts transitively on the
clopens $F(\U_0)$ of $\Q_n$. Together with (1) above
it follows that the diameters and separations of the clopens in $F(\U_0)$
are all bounded multiples of the diameter and separation of a single
clopen in $F(\U_0)$, all lying in some fixed interval $[R^\inv,R]$.
Next, note that   the action of $z^k$ on vertices of $T_p$ induces a
bijection from level $0$ vertices to level $k$ vertices; it follows that
the action of $z^k$ on $\Q_p$ induces a bijection from $\U_0$ to $\U_k$,
and so the action of $z^k$ on
$\Q_n$ induces a bijection from $F(\U_0)$ to $F(\U_k)$. Using
(2) above it follows that the diameters and separations of
all the clopens in
$F(\U_k)$ lie in the interval $\lambda^{-k} [a_0 R^\inv,b_0 R]$. This
establishes Claim~\ref{ClaimDiamsSeps}.

The proofs of (1) and (2) rely on some
general results
about infinite cyclic uniform quasisimilarity actions on metric spaces.
First we have the following result from \cite{FM99}, which
occurs in Proposition~3.3, Step~4:

\specialnumber{26}\proclaim{Lemma}
If $f_n \from X \to X${\rm }, $n \in \Z${\rm ,} is a uniform quasisimilarity action
of the infinite cyclic group $\Z$ on a metric space $X${\rm ,} then
there is a unique number $\lambda \in (0,\infinity)$ with the following
property\/{\rm :}\/ for each $n \in \Z$ the map $f_n$ is $\lambda^n \cdot
[K^\inv,K]$\/{\rm -}\/bilipschitz{\rm ,} where $K$ depends only on the quasisimilarity
constant of the action $(f_n)$.
 \endproclaim

The map $n \mapsto \lambda^n$ is called the {\it stretch homomorphism}
of the infinite cyclic uniform quasisimilarity action $(f_n)$. Note that
$(f_n)$ is uniformly bilipschitz if and only if the stretch homomorphism
is trivial.

Next we have a generalization of \cite[Prop.~3.3]{FM99}, giving a topological characterization of when an
infinite cyclic uniform quasisimilarity action is uniformly bilipschitz.
Given a topological space $X$ and an infinite cyclic action $f_n\from
X\to X$, $n \in \Z$, we say that the action is
{\it locally homothetic at the point $x \in X$} if $f_n(x)=x$ for all
$n$, and in either the positive or negative direction points near
$x$ converge to $x$. That is: for some $\epsilon \in \{-1,+1\}$, and for
some neighborhood $U$ of $x$, we have
$$f_n(y) \to x \quad\hbox{as}\quad n \to \epsilon \cdot \infinity,
\quad\hbox{for all}\quad y \in U.
$$
If the neighborhood $U$ can be taken to be all of $X$ then we say that
$f$ is {\it globally homothetic}.

\specialnumber{27}\proclaim{Lemma}
\label{LemmaUniformBilip}
Let $f_n \from X \to X${\rm ,} $n \in \Z$ be a uniform quasisimilarity action
of the infinite cyclic group $\Z$ on a complete metric space $X$. Then
exactly one of the following happens\/{\rm :}
\begin{itemize}
\item The action $f_n$ is uniformly bilipschitz{\rm ,} with bilipschitz
constant depending only on the quasisimilarity constant.
\item There exists a point $x$ at which the action is locally homothetic.
When this occurs{\rm ,} the local homothety point $x$ is unique{\rm ,} and in fact
$f$ is globally homothetic.
\end{itemize}

\endproclaim
 
{\it Proof}.
If there exists a point of local homothety $x$ for the action $f_n$ then,
for points $y$ near $x$, the ratios

$$\frac{d(f_n(x), f_n(y))}{d(x,y)}, n \in \Z
$$ 
converge to $0$ for $n \to -\infinity$ or $+\infinity$, and so the
ratios
$$\frac{d(f_{-n}(x), f_{-n}(y))}{d(x,y)}
$$
are not bounded, proving that $f_n$ is not uniformly bilipschitz.

Conversely, suppose $f_n$ is not uniformly bilipschitz, with stretch
homomorphism $n\mapsto\lambda^n$. Choose an interval $[a,b]
\subset (0,\infinity)$ so that $f_n$ is $\lambda^n \cdot
[a,b]$-bilipschitz for all $n$. Since $\lambda \ne 1$ it follows that for
sufficiently large $n$ we have $1 \not\in \lambda^n \cdot [a,b]$ and so
there is at most one fixed point. The Contraction Mapping Theorem shows
the existence of a unique fixed point $x$, and it is evident that
$f_n$ is locally homothetic at $x$, indeed $f_n$ is globally homothetic.
\hfill\qed
\vglue12pt

Now we prove (1) and (2) above.

To prove (1), note that the action of $G_0$ on $\Q_p$
is uniformly
bilipschitz, indeed it is isometric. It follows that no cyclic subgroup
of $G_0$ has a local homothety point in $\Q_p$. The actions of $G_0$ on
$\Q_p$ and on $\Q_n$ are topologically conjugate, and so no cyclic
subgroup of $G_0$ has a local homothety point in $\Q_n$. From
Lemma~\ref{LemmaUniformBilip} it follows that each cyclic subgroup of
$G_0$ is uniformly bilipschitz, with bilipschitz constant independent of
the cyclic subgroup since the quasisimilarity constant is independent. In
other words, the action of $G_0$ on $\Q_n$ is uniformly bilipschitz.

To prove (2), let $k \mapsto \lambda^k$ be the stretch
homomorphism of the
action of $Z = \langle z \rangle $ on $\Q_n$. Since $z$
{\it does} have a local homothety point in $\Q_p$ it also has one in
$\Q_n$, and so $\lambda \ne 1$. By replacing $z$ with $z^\inv$ if
necessary we obtain $\lambda>1$.

This completes the proof of Claim~\ref{ClaimDiamsSeps} and of
Theorem~7. The proof of Corollary~8  is
immediate.  \enddemo

As remarked in the introduction, the motivation for
Theorem~7 comes from the end of \cite{FM99}
where it is asked whether a uniform quasisimilarity action on $\Q_n$ is
always conjugate to a similarity action, at least when $n$ is not a
proper power. But this is false. Consider for example $\Sim(\Q_4)$, the
full similarity group of $\Q_4$. Conjugating by a bilipschitz
homeomorphism $\Q_4 \to \Q_2$ we obtain a faithful, uniform
quasisimilarity action of $\Sim(\Q_4)$ on $\Q_2$. However, this action is
not bilipschitz conjugate to a uniform similarity action of $\Sim(\Q_4)$
on
$\Q_2$. To see why, $\Sim(\Q_4)$ contains a copy of the symmetric group
on 4 symbols, a finite group of order $4!=24$. Every \pagebreak similarity action of
a finite group on $\Q_2$ is an isometric action. But every finite subgroup
of the isometry group of $\Q_2$ is a 2~group, because it acts by
direction preserving isometries on the tree $T_2$ fixing some vertex.
\vglue-9pt

\section{Application: Virtually free, cocompact lattices}
\label{SectionFreeLattices}
%\input{QTPartI_Lattices.tex}
% QTPartI_Lattices.tex
\vglue-9pt

In this section we prove Theorem~\ref{TheoremLattices}
which shows that a locally compact group containing a cocompact lattice
which is free of finite rank is closely related to the automorphism group
of a bounded valence, bushy tree, and
Corollary~\ref{CorollaryNoGoodModel} which gives simple examples of
virtually free groups that cannot be cocompact lattices in the same
locally compact group.

We want to understand when two virtually free groups can both act
properly discontinuously and cocompactly by isometries on the same proper
metric space $X$.  Slightly more generally, we ask whether
they are both cocompact lattices in the same locally compact group.

\vglue5pt
{\elevensc Lemma 28.}
{\it Let $\G$ be a locally compact topological group{\rm ,} $\Gamma$ a finitely
generated{\rm ,} cocompact lattice in $\G$, and equip $\Gamma$ with the word
metric. Then $\G$ quasi\/{\rm -}\/acts coboundedly on $\Gamma$.}
\vglue5pt

{\it Proof}.
We prove this by constructing a metric on $\G$ which is left invariant and
so that the inclusion map $\Gamma \inject \G$ is a quasi-isometry,
with respect to the word metric $d_A$ for a fixed finite generating set
$A$ of $\Gamma$. We may then quasiconjugate the left action of $\G$ on
itself to a quasi-action on $\Gamma$ using the inclusion map $\Gamma \to
\G$. The metric we construct on $\G$ does not induce the given topology on
$\G$, but as we are only concerned with the large scale geometry of $\G$
this is not important.

Let $K$ be a symmetric compact set in $\G$ containing the identity and
containing the generating set $A$ for $\Gamma$, so that $\Gamma K =
\G$. It follows that $K$ is a generating set for $\G$. Let $d_K$ be the
left invariant word metric on $\G$ defined in the usual manner using the
generating set $K$: $d_K(g,g')$ is the minimal $n \ge 0$ for which there
exists a $K$-chain of length $n$, meaning a sequence
$g=g_0, g_{1}, \ldots ,g_{n}=g'$ with $g_{i+1}\in g_{i}K$ for all $i$. By
the choice of $K$, every element of $\G$ is within $d_K$ distance $1$ of
an element of $\Gamma$. Thus we need only show that on $\Gamma$ the word
metric $d_A$ and the restriction of $d_{K}$ are quasi-isometric.

Since $K$ contains the generating set $A$ for $\Gamma$, any path in the
Cayley graph with respect to $A$ is a $K$-chain, and  so $d_{K}$ is
bounded above by the word metric~$d_A$.

Given any $K$-chain $\{g_{i}\}$ from $\gamma_{1}$ to $\gamma_{2}$,
define a sequence $\{\gamma_{i}\}$ in $\Gamma$ by choosing, for each $i$,
an  element of $g_{i}K \intersect \Gamma$.  The intersection is nonempty
by the definition of $K$.  Now $\gamma_{i+1} \in g_{i+1}K \subset
g_{i}K^{2} \subset \gamma_{i}K^{3}$, so this sequence in $\Gamma$ is  a
path in the Cayley graph with respect to the generating set $K^{3}
\intersect \Gamma$, which is finite as $K^3$ is compact. Thus $d_{K}$ is
bounded below by the word metric for $K^3 \intersect \Gamma$.

We have shown $d_{K}$ is pinched between two finitely generated word
metrics on $\Gamma$, which proves that it is in the same quasi-isometry
class as the word metrics.
\hfill\qed\vglue4pt

{\it Proof  of Theorem {\rm \ref{TheoremLattices}}}. Given a locally
compact group $\G$ which has a free group $\Gamma$ as a cocompact
lattice, we need to construct a tree $T$ with a $\G$ action, so that the
induced map $\phi \from \G \to \Isom(T)$ is continuous, closed, has
compact kernel, and has cocompact image.

Lemma~28 gives a cobounded quasi-action of $\G$ on the
Cayley graph of~$\Gamma$. Moreover, the proof shows that if we fix a
compactly generated word metric on $\G$ then the left action of $\G$ on
itself is quasiconjugate to resulting quasi-action of $\G$ on the Cayley
graph of $\Gamma$. Theorem~1 produces the tree
$T$ with a $\G$ action. This $\G$ action is quasiconjugate to the
quasi-action of $\G$ on the Cayley graph of $\Gamma$, and so the left
action of $\G$ on itself is quasiconjugate to the action of $\G$ on $T$.
It remains to show using this quasiconjugacy that the map $\phi \from \G
\to \Isom(T)$ has all the desired  properties. For example, the left
action of $\G$ on itself is cobounded, implying that the $\G$-action on
$T$ is cobounded, and so the image of $\phi\from \G
\to\Isom(T)$ is cocompact. 
\vglue6pt
{\elevensc Lemma 29.} 
{\it Given a bounded valence{\rm ,} bushy tree $T${\rm ,} a sequence $(f_{i})$ converges in
$\Isom(T)$ if and only if $(f_i)$ satisfies the following property\/{\rm :}  
 \vglue4pt
{\rm Coarse convergence.} There is some $D$ so that
for any $v$ there is $n$ so that $\{f_{i}(v) \suchthat i \geq  n\}$ has
diameter at most $D$.}
\vglue6pt

{\it Proof}.
``Only if'' is obvious.

Using the local finiteness of $T$, pass to a subsequence of the $f_{i}$
which converges to some $f$.  Further, if, for some $v$, there are
infinitely  many $f_{i}$ for which $f_{i}(v) \neq f(v)$, then we can find
another  convergent    subsequence with a limit, $f'$, different from $f$.
Clearly  one has that $d(f,f') \leq D$.  As isometries of $T$ are unique
in their  bounded distance classes, this is a contradiction proving the
lemma.
\hfill\qed
\vglue6pt

Continuing the proof of Theorem~\ref{TheoremLattices}, consider a
sequence of elements $\{g_{i}\}$ in $\G$. Given two quasiconjugate
quasi-actions of $\G$, the sequence $g_i$ satisfies the coarse convergence
condition of Lemma~29 for one of the quasi-actions if and
only if it satisfies the condition for the other quasi-action. In
particular we can apply this to the left action of $\G$ on itself and to
the quasiconjugate action of $\G$ on $T$. If $g_i$ converges in $\G$
then it satisfies coarse convergence with respect to the left action, and
so $\phi(g_i)$ satisfies coarse convergence in $T$; applying
Lemma~29 it follows that $g_i$ converges in
$\Isom(T)$. This proves that $\phi$ is continuous.

Also, $\phi$ is proper, for take a compact subset $C \subset
\Isom(T)$ and a sequence $g_i$ in $\phi^\inv(C)$. Passing to a
subsequence, $\phi(g_i)$ converges in $C$, and so by
Lemma~29, $\phi(g_i)$ satisfies coarse convergence in $T$.
This implies that $g_i$ satisfies coarse convergence in $\G$. Passing to
a subsequence it follows that $g_i$ converges in $\G$; let $g$ be the
limit. Continuity implies that $\phi(g) = \lim \phi(g_i)$ and so $\phi(g)
\in C$ and $g \in \phi^\inv(C)$.

The kernel of $\phi$, by definition, consists of those elements which
act trivially on $T$. This means that there is some $R$ so that, in the
quasi-action on $\Gamma$, every element of the kernel moves no  point
more than $R$. This means, as in the proof of Lemma~28,
that there is a compact set $K$ in $\G$ so that for  any $g$ in the kernel
of $\phi$, and any $g' \in \G$, $gg' \in g'K$. In particular, $g \in K$.
As the kernel of $\phi$ is closed and contained in a compact set, it is
compact.
 
\vglue6pt
{\it Proof  of Corollary {\rm \ref{CorollaryNoGoodModel}}}.
If two virtually free groups $G$ and $G'$ are both cocompact lattices  in
the same locally compact group, Theorem \ref{TheoremLattices} produces a
tree $T$ on which they both act properly discontinuously  and
cocompactly.  

Consider a group $G={\Bbb  Z}/p{\Bbb  Z} *  {\Bbb 
Z}/p{\Bbb  Z}$ for a prime $p$.  Let $G$ act properly and cocompactly
on a bounded valence, bushy tree $T$, with all vertices of valence at
least three.  We claim $T$ must be the Bass-Serre tree arising from the
given splitting of $G$.

\pagegoal=52pc
Let $X$ be the quotient graph of groups.  We must have $X$ a tree, as
$G$ has no surjections to $\Z$.  Consider an extreme vertex $x$ of $X$,
and let $e$ be the incident edge.  As there are no valence one vertices
in $T$, the edge group of $e$ must be a proper subgroup of the vertex
group of $x$.  Since the action on $T$ is proper, all of the groups are
finite.  Every nontrivial finite subgroup of $G$ is contained in, and hence
equal to, a conjugate of one of the free factors.  Conversely, every such
conjugate must fix a vertex of $T$.  Thus $X$ has exactly two extreme
points, corresponding to the free factors of $G$, with the incident
edges, and all other vertices of $X$, having trivial stabilizers.  A
finite tree with two extreme points is a subdivided interval, and as
$T$ has no valence two vertices, $X$ is a single edge, as claimed.

The groups $\Z/p\Z*\Z/p\Z$ and $\Z/q\Z*\Z/q\Z$ for distinct primes
$p$ and $q$, do not have isomorphic Bass-Serre trees, so the claim
above and Theorem~\ref{TheoremLattices} prove that they are not both
cocompact lattices in any locally compact group.
\hfill\qed

\AuthorRefNames [MSW02b]


\begin{references}
%\pagegoal=48pc

\bibitem{Bie72}
\name{R.~Bieri}, Gruppen mit Poincar\'e-Dualit\"at, {\it Comment.\ Math.\
Helv\/}.\ 
  {\bf 47} (1972), 373--396.

\bibitem{BK90}
\name{H.~Bass} and \name{R.~Kulkarni}, Uniform tree lattices, {\it J.\ Amer.\ Math.\
Soc\/}.\  {\bf 3}
  (1990),  843--902.

\bibitem{BM97}
\name{M.~Burger} and \name{S.~Mozes}, Finitely presented simple groups and products of
  trees, {\it C.\ R.\ Acad.\ Sci.\ Paris S\'er.\ I Math\/}.\
  {\bf 324} (1997),
  747--752.

\bibitem{Bow02}
\name{B.~Bowditch}, Groups acting on Cantor sets and the end structure of
  graphs, {\it Pacific J.\ Math\/}.\  {\bf 207} (2002), 31--60.

\bibitem{BP00}
\name{M.\  Bourdon} and \name{H.\  Pajot}, Rigidity of quasi-isometries for some
  hyperbolic buildings, {\it Comment.\ Math.\ Helv\/}.\ {\bf 75} (2000),
  701--736.

\bibitem{Bro82}
\name{K.~Brown}, {\it Cohomology of Groups\/}, {\it Grad.\ Texts in Math\/}.\ {\bf
87},
Springer-Verlag, New York, 1982.

\bibitem{CC92}
\name{J.~Cannon} and \name{D.~Cooper}, A characterization of cocompact hyperbolic and
  finite-volume hyperbolic groups in dimension three, {\it Trans.\
  Amer.\ Math.\ Soc\/}.\  {\bf 330}
  (1992), 419--431.

\bibitem{Cho96}
\name{R.~Chow}, Groups quasi-isometric to complex hyperbolic space, {\it
Trans.\ Amer.\ Math.\ Soc\/}.\
  {\bf 348} (1996),  1757--1769.

\bibitem{Dav98}
\name{M.~W.\ Davis}, The cohomology of a Coxeter group with group ring
  coefficients, {\it Duke Math.\ J\/}.\  {\bf 91} (1998), 297--314.

\bibitem{Dun85}
\name{M.~J.\ Dunwoody}, The accessibility of finitely presented groups, {\it
Invent.\ Math\/}.\ 
  {\bf 81} (1985), 449--457.

\bibitem{Far97}
\name{B.~Farb}, The quasi-isometry classification of lattices in semisimple
  {Lie} groups, {\it Math.\ Res.\ Lett\/}.\  {\bf 4} (1997),  705--717.

\bibitem{FM98}
\name{B.~Farb} and \name{L.~Mosher}, A rigidity theorem for the solvable
  {Baumslag-Solitar} groups (with an appendix by D.\ Cooper),
  {\it Invent.\ Math\/}.\  {\bf 131} (1998),
  419--451.

\bibitem{FM99}
\bibline, Quasi-isometric rigidity for the solvable
  {Baumslag-Solitar} groups.\  {II}, {\it Invent.\ Math\/}.\
  {\bf 137} (1999),
  613--649.

\bibitem{FM00}
\bibline, On the asymptotic geometry of abelian-by-cyclic
  groups, {\it Acta Math\/}.\ {\bf 184} (2000),  145--202.

\bibitem{FM02}
\bibline, The geometry of surface-by-free groups, {\it
Geom.\ Funct.\ Anal\/}.\
 {\bf 12}  (2002), 915--963.
%%Preprint, \textsc{arXiv:math.GR/0008215}.


\bibitem{FS96}
\name{B.~Farb} and \name{R.~Schwartz}, The large-scale geometry of {Hilbert} modular
  groups, {\it J.\ Differential Geom\/}.\  {\bf 44} (1996),  435--478.

\bibitem{FS87}
\name{M.~Freedman} and \name{R.~Skora}, Strange actions of groups on spheres, {\it J.\
Differential Geom\/}.\
   {\bf 25} (1987), 75--98.

\bibitem{Fur01}
\name{A.~Furman}, Mostow-Margulis rigidity with locally compact targets,
  {\it Geom.\ Funct.\ Anal\/}.\  {\bf 11} (2001),  30--59.

\bibitem{Hin85}
\name{A.~Hinkkanen}, Uniformly quasisymmetric groups, {\it Proc.\ London Math.\
Soc\/}.\ 
  {\bf 51} (1985), 318--338.

\bibitem{JW72}
\name{F.~E.~A. Johnson} and \name{C.~T.~C. Wall}, On groups satisfying Poincar\'e
  duality, {\it Ann.\ of Math\/}.\  {\bf 96} (1972), 592--598.

\bibitem{KK99}
\name{M.~Kapovich} and \name{B.~Kleiner}, Coarse {A}lexander duality and duality
  groups, preprint, arXiv:math.GT/9911003, 1999.

\bibitem{KL97a}
\name{M.~Kapovich} and \name{B.~Leeb}, Quasi-isometries preserve the geometric
  decomposition of {H}aken manifolds, {\it Invent.\ Math\/}.\
  {\bf 128} (1997),
   393--416.

\bibitem{KL97b}
B.~Kleiner and B.~Leeb, Rigidity of quasi-isometries for symmetric spaces
  and Euclidean buildings, {\it IHES Publ.\ Math\/}.\
  {\bf 86} (1997), 115--197.

\bibitem{MSW00}
\name{L.~Mosher, M.~Sageev}, and \name{K.~Whyte}, Quasi-actions on trees, research
  announcement, preprint,  arXiv:math.GR/0005210, 2000.

\bibitem{MSW02a}
\name{L.~Mosher, M.~Sageev}, and \name{K.~Whyte}, Maximally symmetric trees,
  {\it Geom.\ Dedicata\/} {\bf 92} (2002), 195--233.
%%Preprint,   \textsc{arXiv:math.GR/0012004}.

\bibitem{MSW02b}
\bibline, Quasi-actions on trees II:
  Bass-Serre trees, in preparation, 2003.

\bibitem{Pan89a}
\name{P.~Pansu}, Dimension conforme et sph\`ere \`a l'infini des vari\'et\'es \`a
  courbure n\'egative, {\it Ann.\ Acad.\ Sci.\ Fenn.\ Ser.\ A I
  Math\/}.\  {\bf 14}
  (1989),  177--212.

\bibitem{Pan89b}
\bibline, M\'etriques de {Carnot-Carath\'eodory} et quasiisom\'etries des
  espaces\break sym\'etriques de rang un, {\it Ann.\ of Math\/}.\
  {\bf 129} (1989), 1--60.

\bibitem{Pap02}
\name{P.~Papasoglu}, Group splittings and asymptotic topology, preprint,
   arXiv:\break math.GR/0201312, 2002.

\bibitem{Pau95}
\name{F.~Paulin},  De la G\'eom\'etrie et la Dynamique des Groupes Discrets,
  M\'emoire d'Habilitation \`a Diriger les Recherches, ENS Lyon, Juin 1995.

\bibitem{Pau96}
\bibline, Un groupe hyperbolique est d\'etermin\'ee par son bord,
  {\it J.~London Math.\ Soc\/}.\  {\bf 54} (1996),  50--74.

\bibitem{Rei02}
\name{A.~Reiter}, The large scale geometry of products of trees, {\it Geom.\
Dedicata\/} {\bf 92} (2002), 179--184.

\bibitem{Rie01}
\name{E.\ G.\ Rieffel}, Groups quasi-isometric to ${\bf H}\sp 2\times{\bf R}$, {\it J.\  London Math.\ Soc\/}.\  {\bf
64} (2001),  44--60.

\bibitem{Ser80}
\name{J-P.\ Serre}, {\it Trees}, Springer-Verlag, New York, 1980.

\bibitem{Sta68}
\name{J.~Stallings}, On torsion-free groups with infinitely many ends, {\it
Ann.\ of Math\/}.\
   {\bf 88} (1968), 312--334.

\bibitem{Sul}
\name{D.~Sullivan}, private correspondence.

\bibitem{Sul81}
\bibline, On the ergodic theory at infinity of an arbitrary discrete
  group of hyperbolic motions, in {\it Riemann Surfaces and Related
  Topics\/}
  ({\it Proc.\ of the\/} 1978 {\it Stony Brook Conference\/}),
  {\it Ann.\ of Math.\
  Studies\/} {\bf 97},
  Princeton Univ.\ Press, Princeton, NJ, 1981, 465--496.

\bibitem{SW79}
\name{P.~Scott} and \name{C.~T.~C. Wall}, Topological methods in group theory, in
  {\it Homological Group Theory\/} ({\it Proc.\ of Durham Sympos\/}.,
  {\it Sept\/}.\
  1977), {\it London
  Math.\ Soc.\ Lecture Notes\/} {\bf 36} (1979), 137--203, Cambridge
  Univ.\ Press, New York.

\bibitem{Tab00}
\name{J.~Taback}, Quasi-isometric rigidity for ${\rm PSL}(2,{\bf Z}[1/p])$, {\it Duke Math.\ J\/}.\  {\bf 101} (2000),  335--357.

\bibitem{Tuk80}
\name{P.~Tukia}, On two-dimensional quasiconformal groups, {\it Ann.\ Acad.\
Sci.\ Fenn.\ Ser.\ A I Math\/}.\
 {\bf 5} (1980), 73--78.

\bibitem{Tuk81}
\bibline, A quasiconformal group not isomorphic to a M\"obius group,
  {\it Ann.\ Acad.\ Sci.\ Fenn.\ Ser.\ A I Math\/}.\
  {\bf 6} (1981), 149--160.


\bibitem{Tuk86}
\bibline, On quasi-conformal groups, {\it J.\ Analyse Math\/}.\  {\bf 46}
  (1986), 318--346.

\bibitem{Why}
\name{K.~Whyte}, Geometries which fiber over trees, in preparation.

\bibitem{Why02}
\bibline, The large scale geometry of the higher Baumslag-Solitar
  groups, {\it Geom.\ Funct.\ Anal\/}.\  {\bf 11} (2002), 1327--1353.

\end{references}
\end{document}